\documentclass{article}
\usepackage{amssymb}
\usepackage{graphicx}
\usepackage{amsmath}

\input{tcilatex}
\begin{document}

\begin{center}
\bigskip

\bigskip \textbf{VELI SHAKHMUROV}

\bigskip

\textbf{EMBEDDING OPERATORS IN\ SOBOLEV-LIONS\ SPACES AND\ APPLICATIONS}

\bigskip

Department of Mechanical Engineering, Okan University, Akfirat, Tuzla 34959
Istanbul, Turkey,

\textbf{AMS: \ 26Bxx, 43Axx, 46Bxx, 35Jxx}

\bigskip \textbf{Abstract}
\end{center}

The embedding theorems in Sobolev-Lions type anisotropic weighted spaces $W_{%
\mathbf{p},\gamma }^{l}\left( \Omega ;E_{0},E\right) $ are studied, here $%
E_{0}$ and $E$ are two Banach spaces. The most regular interpolation spaces $%
\left( E_{0},E\right) _{p,\theta _{\alpha }}$ between $E_{0}$ and $E$ are
found such that the mixed differential operators $D^{\alpha }$ are bounded
from $W_{\mathbf{p},\gamma }^{l}\left( \Omega ;E_{0},E\right) $ to $L_{%
\mathbf{p},\gamma }\left( \Omega ;\left( E_{0},E\right) _{p,\theta _{\alpha
}}\right) ,$ where $L_{\mathbf{p},\gamma }$ denotes weighted abstact
Lebesgue space\ with mixed nom and%
\begin{equation*}
\alpha =\left( \alpha _{1},\alpha _{2},...,\alpha _{n}\right) ,l=\left(
l_{1},l_{2},...,l_{n}\right) \text{, }\theta _{\alpha
}=\dsum\limits_{k=1}^{n}\frac{\alpha _{k}}{l_{k}}.
\end{equation*}%
By applying this result separability properties of degenerate anisotropic
differential operator equations, well-posedeness and Strichartz type
estimates for solution of corresponding parabolic problem are established.


\begin{center}
\bigskip
\end{center}

\bigskip \textbf{1. Introduction}

Embedding of function spaces were studied in a series of books and papers
(see, for example $\left[ \text{2, 4, 14, 15}\right] $). The embedding
properties of abstract function spaces have been considered e.g. in $\left[ 
\text{1, 6, 12, 16-21}\right] .$ Lions-Peetre $\left[ 12\right] $ showed
that if $\ u\in $ $L_{2}\left( 0,T,H_{0}\right) $ and $u^{(m)}\in
L_{2}(0,T,H $ $),$ then 
\begin{equation*}
u^{(i)}\in L_{2}(0,T,\left[ H,H_{0}\right] _{\frac{i}{m}}),\ i=1,2,...m-1,
\end{equation*}%
where $H_{0}$, $H$ are Hilbert spaces, $H_{0}$ is continuously and densely
embedded into $H$, and $\left[ H_{0},H\right] _{\theta }$\ is an
interpolation space between $H_{0}$, $H$ for $0\leq \theta \leq 1.$ The
similar questions for $H$-valued Sobolev spaces $W_{2}^{l}\left( \Omega
;H_{0},H\right) $ studied in $\left[ 24\right] ,$ where $\Omega \subset
R^{n}.$ Then, the boundedness of differential operator $u\rightarrow
D^{\alpha }u$ from $W_{p}^{l}\left( \Omega ;H_{0},H\right) $ to $L_{p}\left(
\Omega ;\left( H_{0},H\right) _{p,\left\vert \alpha :l\right\vert }\right) $
were considered in $\left[ \text{18, 19}\right] .$ This question is
generalized for corresponding weighted spaces in $\left[ 13\right] $. Later,
such type embedding results in $E$-valued function spaces $W_{p,\gamma
}^{l}\left( \Omega ;E_{0},E\right) $ and its weighted versions studied in $%
\left[ 18-21\right] $. In this paper, we prove the continuity and
compactness\ of embedding operators in weighted anisotropic function spaces\ 
$W_{\mathbf{p},\gamma }^{l}\left( \Omega ;E_{0},E\right) $ for mixed $%
\mathbf{p}$, which will be defined in bellow.

Here $l=\left( l_{1},l_{2},...l_{n}\right) $, $l_{k}$ are positive integers
and $\gamma \left( x\right) $ is a positive measurable function on $\Omega
\subset R^{n}$. Let $\alpha _{1},\alpha _{2},...,\alpha _{n}$ be nonnegative
integers, $\mathbf{p}=\left( p_{1},p_{2},...,p_{n}\right) $ and $\mathbf{q}%
=\left( q_{1},q_{2}...,q_{n}\right) $ for $1\leq p_{k}\leq q_{k}<\infty ,$%
\begin{equation*}
\varkappa =\sum\limits_{k=1}^{n}\frac{\alpha _{k}+\frac{1}{p_{k}}-\frac{1}{%
q_{k}}}{l_{k}},\text{ }D^{\alpha }=D_{1}^{\alpha _{1}}D_{2}^{\alpha
_{2}}...D_{n}^{\alpha _{n}}=\frac{\partial ^{\alpha }}{\partial
x_{1}^{\alpha _{1}}\partial x_{2}^{\alpha _{2}}...\partial x_{n}^{\alpha
_{n}}}.
\end{equation*}

Let\ $A$\ be a positive operator in $E$,\ then there are fractional powers
of the operator $A\ $(see $\left[ 22\right] ,$ \S 1.15.1) and for each
powers $A^{\theta }$ of $A$ let $E\left( A^{\theta }\right) $ denote the
domain $D\left( A^{\theta }\right) $ of $A^{\theta }$ with graphical norm.
Under certain assumptions to be stated later, we prove that differential
operators $u\rightarrow D^{\alpha }u$ are bounded from $W_{\mathbf{p},\gamma
}^{l}\left( \Omega ;E(A\right) ,E)$ to $L_{\mathbf{q},\gamma }\Omega ;\left(
E\left( A\right) ,E\right) _{\varkappa +\mu ,\sigma }),$ i.e embedding \ \ 

\begin{equation*}
D^{\alpha }W_{\mathbf{p},\gamma }^{l}\left( \Omega ;E\left( A\right)
,E\right) \subset L_{\mathbf{q},\gamma }\left( \Omega ;\left( E\left(
A\right) ,E\right) _{\varkappa +\mu ,\sigma }\right)
\end{equation*}%
is continuous. More precisely, we prove the following uniform sharp estimate

\begin{equation*}
\left\Vert D^{\alpha }u\right\Vert _{L_{\mathbf{q},\gamma }\Omega ;\left(
E\left( A\right) ,E\right) _{\varkappa +\mu ,\sigma })}\leq C_{\mu }(h^{\mu
}\left\Vert u\right\Vert _{W_{\mathbf{p},\gamma }^{l}\left( \Omega
;.E(A)E\right) }+h^{-(1-\mu )}\left\Vert u\right\Vert _{L_{\mathbf{p},\gamma
}(\Omega ,;E)})
\end{equation*}%
for $u\in W_{\mathbf{p},\gamma }^{l}\left( \Omega ;E(A),E\right) $, $0\leq
\mu \leq 1-\varkappa $ and $h>0$, where $\sigma =\max\limits_{k}p_{k}.$ The
constant $C_{\mu }$ is independent of $u$ and$\ h.$ Further, we prove
compactness of embedding operator. These kind of embedding theorems occur in
the investigation of boundary value problems for anisotropic elliptic
differential-operator equations%
\begin{equation}
\ \sum\limits_{k=1}^{n}(-1)^{l_{k}}t_{k}\frac{\partial ^{2l_{k}}u}{\partial
x_{k}^{2l_{k}}}+\left( A+\lambda \right) u\left( x\right)
+\sum\limits_{\left\vert \alpha :2l\right\vert
<1}\prod\limits_{k=1}^{n}t_{k}^{\frac{\alpha _{k}}{2l_{k}}}A_{\alpha
}(x)D^{\alpha }u\left( x\right) =f\left( x\right) ,  \tag{1.1}
\end{equation}%
$\ $where\ $A$, $A_{\alpha }\left( x\right) $ are linear operators in a
Banach space $E$, $\ t_{k}$ are positive parameters and $\lambda $ is a
complex number. By using the above embedding results and operator valued
multiplier theorems we obtain that problem $\left( 1.1\right) $ is uniform
separable in $L_{\mathbf{p,\gamma }}\left( R^{n};E\right) ,$ i.e. for $f\in $
$L_{\mathbf{p,\gamma }}\left( R^{n};E\right) $ problem $\left( 1.1\right) $
has a unique solution $u\in W_{\mathbf{p,\gamma }}^{2l}(R^{n};(E(A),E)$ and
the following uniform coercive estimate holds%
\begin{equation*}
\ \sum\limits_{k=1}^{n}t_{k}\left\Vert \frac{\partial ^{2l_{k}}u}{\partial
x_{k}^{2l_{k}}}\right\Vert _{L_{\mathbf{p,\gamma }}\left( R^{n};E\right)
}+\left\Vert Au\right\Vert _{L_{\mathbf{p,\gamma }}\left( R^{n};E\right)
}\leq C\left\Vert f\right\Vert _{L_{\mathbf{p,\gamma }}\left( R^{n};E\right)
},
\end{equation*}%
where the constant $C$ depend only on $\mathbf{p}$ and $l.$

Moreover, we get the following uniform sharp resolvent estimate 
\begin{equation*}
\sum\limits_{k=1}^{n}\sum\limits_{i=0}^{2l_{k}}t_{k}^{\frac{i}{2l_{k}}%
}\left\vert \lambda \right\vert ^{1-\frac{i}{2l_{k}}}\left\Vert \frac{%
\partial ^{i}}{\partial x_{k}^{i}}\left( O_{t}+\lambda \right)
^{-1}\right\Vert +\left\Vert A\left( O_{t}+\lambda \right) ^{-1}\right\Vert
\leq C,
\end{equation*}%
where $O_{t}$ is the operator generated by problem $\left( 1.1\right) .$

For $l_{1}=l_{2}=,...,=l_{n}=m$ we get the elliptic differential-operator
equation%
\begin{equation*}
\ \sum\limits_{k=1}^{n}(-1)^{m}t_{k}\frac{\partial ^{2m}u}{\partial
x_{k}^{2m}}+Au\left( x\right) +\sum\limits_{\left\vert \alpha \right\vert
<2m}\prod\limits_{k=1}^{n}t_{k}^{\frac{\alpha _{k}}{2m}}A_{\alpha
}(x)D^{\alpha }u\left( x\right) =f\left( x\right) .
\end{equation*}

Then, by using regularity properties of $\left( 1.1\right) $ the
well-posedeness and uniform Strichartz type estimates are established for
the solution of abstract parabolic problem 
\begin{equation}
\frac{\partial u}{\partial t}\ +\sum\limits_{k=1}^{n}(-1)^{l_{k}}\varepsilon
_{k}\frac{\partial ^{2l_{k}}u}{\partial x_{k}^{2l_{k}}}+Au=f\left(
t,x\right) ,\text{ }  \tag{1.2}
\end{equation}

\begin{equation*}
u\left( 0,x\right) =0\text{, }x\in R^{n},\text{ }t\in \left( 0,T\right) ,
\end{equation*}
where $A$ is a linear operator in a Banach space $E$ and $\varepsilon _{k}$
are small positive parameters.

In this direction we should mention e.g.\ the works presented in $\left[
18-21\right] ,$ $\left[ 22\right] $, $\left[ 26\right] $.

Modern analysis methods, particularly abstract harmonic analysis, the
operator theory, interpolation of Banach spaces, theory of semigroups and
perturbation theory of linear operators are the main tools implemented to
carry out the analysis.

\begin{center}
\ \ \textbf{2. Notations and definitions}
\end{center}

Let $\mathbb{R}$, $\mathbb{C}$ be the sets of real and complex numbers,
respectively. Let $E_{1}$ and $E_{2}$ be two Banach spaces and $L\left(
E_{1},E_{2}\right) $ denotes the spaces of bounded linear operators from $%
E_{1}$ to $E_{2}.$ For $E_{1}=E_{2}=E$ we denote $L\left( E,E\right) $ by $%
L\left( E\right) .$ We will sometimes write $A+\xi $ instead of\ $A+\xi I$\
for a scalar $\xi $ and $\left( A+\xi I\right) ^{-1}$ denotes the resolvent
of operator $A$, where $I$ is the identity operator in $E.$

Let 
\begin{equation*}
S_{\varphi }=\left\{ \text{ }\xi \in \mathbb{C},\left\vert \text{arg }\xi
\right\vert \leq \varphi \right\} \cup \left\{ 0\right\} ,\text{ }0<\varphi
\leq \pi .
\end{equation*}

\textbf{Definition 1}. A linear operator\ $A$\ is said to be $\varphi -$%
positive\ in a Banach space\ $E,$\ if \ $D\left( A\right) $ is dense on $E$
and

\begin{equation*}
\left\Vert \left( A+\xi I\right) ^{-1}\right\Vert _{L\left( E\right) }\leq
M\left( 1+\left\vert \xi \right\vert \right) ^{-1}
\end{equation*}

with\ $\xi \in S_{\varphi }$, where\ $M$ is a positive constant.\ 

\bigskip \textbf{Definition 2. }For $-\infty <\theta <\infty $ let%
\begin{equation*}
E\left( A^{\theta }\right) =\left\{ u\in D\left( A^{\theta }\right)
,\left\Vert u\right\Vert _{E\left( A^{\theta }\right) }=\left\Vert A^{\theta
}u\right\Vert _{E}+\left\Vert u\right\Vert _{E}<\infty \right\} .\ \ 
\end{equation*}

Let $\gamma =\gamma \left( x\right) $ be a positive measurable function on $%
\Omega \subset R^{n}.$

\textbf{Definition 3}. $L_{p,\gamma }\left( \Omega ;E\right) $ denotes the
space of strongly measurable $E$-valued functions such that are defined on $%
\Omega $ with the norm

\begin{equation*}
\left\Vert u\right\Vert _{L_{p,\gamma }\left( \Omega ;E\right) }=\left(
\int_{\Omega }\left\Vert u\right\Vert _{E}^{p}\gamma \left( x\right)
dx\right) ^{1/p}<\infty ,\text{ }1\leq p<\infty .
\end{equation*}

For $\gamma \left( x\right) \equiv 1$ we denote $L_{p,\gamma }\left( \Omega
;E\right) $ by $L_{p}\left( \Omega ;E\right) .$\ For $\mathbf{p}=\left(
p_{1},p_{2},...,p_{n}\right) ,$ $1\leq p_{k}<\infty $ we denote by $L_{%
\mathbf{p,\gamma }}\left( R^{n};E\right) $ the space of all $E$-valued
strongly measurable on $R^{n}$ functions with mixed norm%
\begin{equation*}
\left\Vert f\right\Vert _{L_{\mathbf{p,\gamma }}\left( R^{n},E\right) }=
\end{equation*}

\begin{equation*}
\left( \left( \int\limits_{-\infty }^{\infty }\left( ...\int\limits_{-\infty
}^{\infty }\left( \int\limits_{-\infty }^{\infty }\left\Vert f\left(
x\right) \right\Vert _{E}^{p_{1}}\gamma \left( x\right) dx_{1}\right) ^{%
\frac{p_{2}}{p_{1}}}dx_{2}\right) ^{\frac{p_{3}}{p_{2}}}...\right) ^{\frac{%
p_{n}}{p_{n-1}}}dx_{n}\right) ^{\frac{1}{p_{n}}}<\infty .
\end{equation*}%
For $\gamma \left( x\right) \equiv 1$ we denote $L_{\mathbf{p},\gamma
}\left( R^{n};E\right) $ by $L_{\mathbf{p}}\left( R^{n};E\right) .$

The weight $\gamma \left( x\right) $ satisfies an $A_{p}$ condition; i.e., $%
\gamma \left( x\right) \in A_{p},$ $p\in \left( 1,\infty \right) $ if there
is a positive constant $C$ such that

\begin{equation*}
\left( \frac{1}{\left\vert Q\right\vert }\dint\limits_{Q}\gamma (x)dx\right)
\left( \frac{1}{\left\vert Q\right\vert }\dint\limits_{Q}\gamma ^{-\dfrac{1}{%
p-1}}(x)dx\right) ^{p-1}\leq C
\end{equation*}%
for all compacts $Q\subset R^{n}$.

\textbf{Remark 2.1. }The result $\left[ 23\right] $ implies that the space $%
l_{p}$, $p\in \left( 1,\infty \right) $ satisfies the multiplier condition
with respect to $p$ and the weight functions 
\begin{equation*}
\gamma =\left\vert x\right\vert ^{\alpha }\text{, }-1<\alpha <p-1\text{, }%
\gamma =\prod\limits_{k=1}^{N}\left( 1+\sum\limits_{j=1}^{n}\left\vert
x_{j}\right\vert ^{\alpha _{jk}}\right) ^{\beta _{k}}\text{, }
\end{equation*}

\begin{equation*}
\alpha _{jk}\geq 0\text{, }N\in \mathbb{N}\text{, }\beta _{k}\in \mathbb{R},
\end{equation*}%
Suppose that\ $S=S\left( R^{n}\right) $ is Schwartzs space of test functions
and\ $S^{^{\prime }}\left( E\right) =S\left( R^{n};E\right) $\ is the space
of linear continued mapping from $S$ into\ $E$ and is called \ $E-$ valued
Schwartzs distributions. For $\varphi \in S$ the Fourier transform $\hat{%
\varphi}$ and inverse Fourier transform $\check{\varphi}$ are defined by the
relations%
\begin{equation*}
\hat{\varphi}\left( \xi \right) =\left( F\varphi \right) \left( \xi \right)
=\left( 2\pi \right) ^{-n/2}\dint\limits_{R^{n}}\varphi \left( x\right)
e^{-i\xi x}dx,
\end{equation*}

\begin{equation*}
\check{\varphi}\left( x\right) =\left( F^{-\shortmid }\varphi \right) \left(
x\right) =\left( 2\pi \right) ^{-n/2}\dint\limits_{R^{n}}\varphi \left( \xi
\right) e^{i\xi x}d\xi ,
\end{equation*}%
where%
\begin{equation*}
\ \xi =(\ \xi _{1},\xi _{2},...,\xi _{n}),\text{ }x=\left(
x_{1},x_{2},...,x_{n}\right) ,\ \ \xi x=\xi _{1}x_{1}+\xi _{2}x_{2}+...+\xi
_{n}x_{n}.
\end{equation*}

The Fourier transformation and the inverse Fourier transformation of \ $E-$%
valued generalized functions \ $f\in S^{^{\prime }}\left( R^{n};E\right) $
are defined by the relations.

\begin{equation*}
<f^{\symbol{94}},\varphi >=<f,\hat{\varphi}>\text{and}<\check{f},\varphi
>=<f,\check{\varphi}>,\ \ 
\end{equation*}%
where \ $<f,\varphi >$ \ means the value of generalized function \ $f\in
S^{^{\prime }}\left( R^{n};E\right) $ on the $\varphi \in S\left(
R^{n}\right) $.

\textbf{Definition 4}. Let\ $\alpha =\left( \alpha _{1},\alpha
_{2},...,\alpha _{n}\right) $, where $\alpha _{i}$ are positive integers.
The $E$ $-$values generalized functions\ $D^{\alpha }f$ is called the
generalized derivatives in the sense of Schwarts distributions of the
generalized function \ $f\in S^{^{\prime }}\left( R^{n};E\right) $ if the
relation

\begin{equation*}
<D^{\alpha }f,\varphi >=\left( -1\right) ^{\left\vert \alpha \right\vert
}<f,D^{r}\varphi >
\end{equation*}%
holds for all $\varphi \in S.$

\bigskip\ It is known for all\ $\varphi \in S$\ the relations \ \ 

\begin{equation}
F\left( D_{x}^{\alpha }\varphi \right) =\left( i\xi _{1}\right) ^{\alpha
_{1}}...\left( i\xi _{n}\right) ^{\alpha _{n}}\hat{\varphi},\text{ }D_{\xi
}^{\alpha }F\left[ \varphi \right] =F\left[ \left( -ix_{1}\right) ^{\alpha
_{1}}...\left( -ix_{n}\right) ^{\alpha _{n}}\varphi \right] \   \tag{2.1}
\end{equation}%
$\ $holds.

\ Let\ $\eta $ be a infinitely differentiable function with polynomial
structure and \ $f\in S^{^{\prime }}\left( R^{n};E\right) $. Then\ $\eta
f\in S^{^{\prime }}\left( R^{n},E\right) $ is a generalized function defined
by the relation

\begin{equation*}
<\eta f,\varphi >=<f,\eta f>\ \ \forall \varphi \in S\left( R^{n}\right) .
\end{equation*}%
\ By using Definition 4 and relations $\left( 2.1\right) $ we get

\begin{equation}
\ F\left( D_{x}^{\alpha }f\right) =\left( i\xi _{1}\right) ^{\alpha
_{1}}...\left( i\xi _{n}\right) ^{\alpha _{n}}f^{\symbol{94}},\ D_{\xi
}^{\alpha }\left( F\left( f\right) \right) =F\left[ \left( -ix_{n}\right)
^{\alpha _{1}}...\left( -ix_{n}\right) ^{\alpha _{n}}f\ \ \right]  \tag{2.2}
\end{equation}%
for \ $f\in S^{^{\prime }}\left( R^{n};E\right) $.

\bigskip The Banach space\ $E$ is called an UMD-space if\ the Hilbert
operator%
\begin{equation*}
\left( Hf\right) \left( x\right) =\lim\limits_{\varepsilon \rightarrow
0}\int\limits_{\left\vert x-y\right\vert >\varepsilon }\frac{f\left(
y\right) }{x-y}dy
\end{equation*}%
\ is bounded in $L^{p}\left( \mathbb{R},E\right) $ for $p\in \left( 1,\infty
\right) $ (see. e.g. $\left[ 2\right] $). UMD spaces include e.g. $L^{p}$, $%
l_{p}$ spaces and Lorentz spaces $L_{pq}$ for $p$, $q\in \left( 1,\infty
\right) $.

$C^{\left( m\right) }\left( \Omega ;E\right) $\ will denote the spaces of $%
E- $valued bounded uniformly strongly $m$ times continuously differentiable
functions on $\Omega .$ Assume $\gamma $ is such that $S\left(
R^{n};E_{1}\right) $ is dense in $L_{\mathbf{p},\gamma }\left(
R^{n};E_{1}\right) .$

\textbf{Definition 5}. A function $\Psi \in C^{\left( m\right) }\left(
R^{n};L\left( E_{1},E_{2}\right) \right) $\ is called a multiplier from\ $L_{%
\mathbf{p},\gamma }\left( R^{n};E_{1}\right) $ to $L_{\mathbf{q},\gamma
}\left( R^{n};E_{2}\right) $\ if the map $u\rightarrow \Lambda u=F^{-1}\Psi
\left( \xi \right) Fu,$ $u\in S\left( R^{n};E_{1}\right) $ is well defined
and extends to a bounded linear operator 
\begin{equation*}
\Lambda :\ L_{\mathbf{p},\gamma }\left( R^{n};E_{1}\right) \rightarrow L_{%
\mathbf{q},\gamma }\left( R^{n};E_{2}\right) .
\end{equation*}

We denote the set of all multipliers from $L_{\mathbf{p},\gamma }\left(
R^{n};E_{1}\right) $ to $L_{\mathbf{q},\gamma }\left( R^{n};E_{2}\right) $
by $M_{\mathbf{p},\gamma }^{\mathbf{q},\gamma }\left( E_{1},E_{2}\right) .$
For $E_{1}=E_{1}=E$ it denotes by\ $M_{\mathbf{p},\gamma }^{\mathbf{q}%
,\gamma }\left( E\right) .$

Let $\Phi _{h}=\left\{ \Psi _{h}\in M_{\mathbf{p},\gamma }^{\mathbf{q}%
,\gamma }\left( E_{1},E_{2}\right) ,\text{ }h\in Q\right\} $ denote a
collection of multipliers depending on the parameter $h.$

We say that $W_{h}$ is a uniform collection of multipliers if there exists a
positive constant $M$ independent of $h\in Q$ such that

\begin{equation*}
\left\Vert F^{-1}\Psi _{h}Fu\right\Vert _{L_{\mathbf{q},\gamma }\left(
R^{n};E_{2}\right) }\leq M\left\Vert u\right\Vert _{L_{\mathbf{p},\gamma
}\left( R^{n};E_{1}\right) }\ \ \ \ \ \ 
\end{equation*}%
for all $u\in S\left( R^{n};E_{1}\right) $ and $h\in Q.$

Note that, Fourier multiplier theorems in complex valued weighted $L_{p}$
spaces investigated e.g. in $\left[ \text{9-11}\right] .$ In Banach
space-valued classes this question studied e.g. in $\left[ \text{5, 6}\right]
.$

Let 
\begin{equation*}
U_{n}=\beta =\left( \beta _{1},\beta _{2},...,\beta _{n}\right) ,\text{ }%
\beta _{i}\in \left\{ 0,1\right\} ,\text{ }\sigma _{k}=\text{ }\left( \frac{1%
}{p_{k}}-\frac{1}{q_{k}}\right) ,
\end{equation*}%
\begin{equation*}
V_{n}=\left\{ \xi =\left( \xi _{1},\xi _{2},...,\xi _{n}\right) \in R^{n},%
\text{ }\xi _{i}\neq 0,\text{ }i=1,2,...,n\right\} .
\end{equation*}%
\textbf{Definition 6. }The Banach space $E$ satisfies the multiplier
condition with respect to\ $\mathbf{p}$ and\ $\mathbf{q}$ (or with respect
to $\mathbf{p}$\ in the case of$\ \mathbf{p}=\mathbf{q}$) and with respect
to weighted function $\gamma $\ if for all $\Psi \in C^{\left( n\right)
}\left( R^{n};L\left( E\right) \right) $ with $\beta \in U_{n},$\ $\xi \in
V_{n}$ the inequality

\begin{equation}
\left\Vert D^{\beta }\Psi \left( \xi \right) \right\Vert _{L\left( E\right)
}\leq C\dprod\limits_{k=1}^{n}\left\vert \xi _{k}\right\vert ^{-\left( \beta
_{k}+\sigma _{k}\right) }  \tag{2.3}
\end{equation}%
implies $\Psi \in M_{\mathbf{p},\gamma }^{\mathbf{q},\gamma }\left( E\right)
.$

Note that, if $E_{1}$ and $E_{2}$ are UMD spaces, $\gamma \left( x\right)
\equiv 1$ and $p_{1}=p_{2}=...=p_{n}=p,$ then by virtue of operator valued
multiplier theorems (see e.g $\left[ 7\right] ,$ $\left[ 24\right] $) we
obtain that $\Psi $ is a Fourier multiplier in $L_{p}\left( R^{n};E\right) .$
It is well known (see $\left[ 12\right] $) that any Hilbert space satisfies
the multiplier condition for $\gamma \left( x\right) \equiv 1$ with respect
to any\ $p$\ and\ $q$ with \ $1<p\leq q<\infty .$ However, there are Banach
spaces which are not Hilbert spaces but satisfy the multiplier condition,
for example UMD spaces, $\xi -$convex Banach lattice spaces (see $\left[ 7%
\right] $, $\left[ 24\right] ,$ $\left[ 5\right] $, $\left[ 6\right] $).

Assume $\gamma _{k}$ are positive measurable functions on $\mathbb{R}$ and 
\begin{equation*}
D_{k}^{\left[ i\right] }=\left( \gamma _{k}\left( x_{k}\right) \frac{%
\partial }{\partial x_{k}}\right) ^{i}.
\end{equation*}

\textbf{Definition 7}. Consider the following spaces: 
\begin{equation*}
W_{\mathbf{p,\gamma }}^{l}\left( \Omega ,E_{0},E\right) =\left\{ u\in L_{%
\mathbf{p},\gamma }\left( \Omega ;E_{0}\right) ,\text{ }D_{k}^{l_{k}}u\in L_{%
\mathbf{p},\gamma }\left( \Omega ;E\right) \text{, }k=1,2,...,n\right\} ,
\end{equation*}%
\begin{equation*}
\left\Vert u\right\Vert _{W_{\mathbf{p,\gamma }}^{l}\left( \Omega
;E_{0,}E\right) }=\left\Vert u\right\Vert _{L_{\mathbf{p,\gamma }}\left(
\Omega ;E_{0}\right) }+\sum\limits_{k=1}^{n}\left\Vert
D_{k}^{l_{k}}u\right\Vert _{L_{\mathbf{p,\gamma }}\left( \Omega ;E\right)
}<\infty ,
\end{equation*}

\begin{equation*}
W_{\mathbf{p,\gamma }}^{\left[ l\right] }\left( \Omega ,E_{0},E\right)
=\left\{ u\in L_{\mathbf{p}}\left( \Omega ;E_{0}\right) ,\text{ }D_{k}^{%
\left[ l_{k}\right] }u\in L_{\mathbf{p}}\left( \Omega ;E\right) \right\} ,
\end{equation*}%
\begin{equation*}
\left\Vert u\right\Vert _{W_{\mathbf{p,\gamma }}^{\left[ l\right] }\left(
\Omega ;E_{0,}E\right) }=\left\Vert u\right\Vert _{L_{\mathbf{p}}\left(
\Omega ;E_{0}\right) }+\sum\limits_{k=1}^{n}\left\Vert D_{k}^{\left[ l_{k}%
\right] }u\right\Vert _{L_{\mathbf{p}}\left( \Omega ;E\right) }<\infty .
\end{equation*}

For $\gamma \left( x\right) \equiv 1$ we will denote $W_{\mathbf{p,\gamma }%
}^{l}\left( \Omega ,E_{0},E\right) $ and $W_{\mathbf{p,\gamma }}^{\left[ l%
\right] }\left( \Omega ,E_{0},E\right) $ by $W_{\mathbf{p}}^{l}\left( \Omega
,E_{0},E\right) $.

Let$\ t_{k}$ be positive parameters and $t=\left(
t_{1},t_{2},...,t_{n}\right) $. We define the following parametrized norms
in $W_{\mathbf{p},\gamma }^{l}\left( \Omega ;E_{0},E\right) $ and $W_{%
\mathbf{p,\gamma }}^{\left[ l\right] }\left( \Omega ,E_{1},E\right) $ such
that 
\begin{equation*}
\left\Vert u\right\Vert _{W_{\mathbf{p,\gamma }}^{l},_{t}\left( \Omega
;E_{0},E\right) }=\left\Vert u\right\Vert _{L_{\mathbf{p,\gamma }}\left(
\Omega ;E_{0}\right) }+\sum\limits_{k=1}^{n}\left\Vert
t_{k}D_{k}^{l_{k}}u\right\Vert _{L_{\mathbf{p,\gamma }}\left( \Omega
;E\right) }<\infty ,
\end{equation*}%
\begin{equation*}
\left\Vert u\right\Vert _{W_{\mathbf{p,\gamma }}^{\left[ l\right]
},_{t}\left( \Omega ;E_{0},E\right) }=\left\Vert u\right\Vert _{L_{\mathbf{p}%
}\left( \Omega ;E_{0}\right) }+\sum\limits_{k=1}^{n}\left\Vert t_{k}D_{k}^{%
\left[ l_{k}\right] }u\right\Vert _{L_{\mathbf{p}}\left( \Omega ;E\right)
}<\infty .
\end{equation*}

For two elements $u$, $\upsilon \in E$ the expression $\left\Vert
u\right\Vert \sim \left\Vert \upsilon \right\Vert $ means that there exist
positive numbers $C_{1}$ and $C_{2}$ such that 
\begin{equation*}
C_{1}\left\Vert u\right\Vert \leq \left\Vert \upsilon \right\Vert \leq
C_{2}\left\Vert u\right\Vert .
\end{equation*}

\begin{center}
\ \textbf{3. Embedding theorems}
\end{center}

In this section, we prove that the generalized derivative operator\ $%
D^{\alpha }$ generates a continuous embedding in Sobolev spaces of
vector-functions. Let $\alpha _{k}$ be nonnegative and $l_{k}$ positive
integers and 
\begin{equation*}
\ \alpha =\left( \alpha _{1},\alpha _{2},...,\alpha _{n}\right) ,\text{ }%
l=\left( l_{1},l_{2},...,l_{n}\right) ,\text{ }\sigma =\left( \sigma
_{1},\sigma _{2},...,\sigma _{n}\right) ,\text{ }
\end{equation*}%
\begin{equation*}
t=\left( t_{1},t_{2},...,t_{n}\right) \text{, }t_{k}>0,\text{ }\sigma _{k}>0,%
\text{ }
\end{equation*}%
\begin{equation*}
1<p_{k}\leq q_{k}<\infty \text{, \ }\mathbf{p}=\left(
p_{1},p_{2},...,p_{n}\right) ,\text{ }\mathbf{q}=\left(
q_{1},q_{2},...,q_{n}\right) ,
\end{equation*}%
\begin{equation*}
\varkappa =\left\vert \left( \alpha +\frac{1}{\mathbf{p}}-\frac{1}{\mathbf{q}%
}\right) :l\right\vert =\dsum\limits_{k=1}^{n}\frac{\alpha _{k}+\sigma _{k}}{%
l_{k}},
\end{equation*}%
\begin{equation*}
\xi =\left( \xi _{1},\xi _{2},...,\xi _{n}\right) \in R^{n},\text{ }%
\left\vert \xi \right\vert ^{\alpha +\sigma
}=\prod\limits_{k=1}^{n}\left\vert \xi _{k}\right\vert ^{\alpha _{k}+\sigma
_{k}},\text{ }T\left( t\right) =\text{ }\prod\limits_{k=1}^{n}t_{k}^{\frac{%
\alpha _{k}+\sigma _{k}}{l_{k}}}.
\end{equation*}%
First of all, in a similar way as in $\left[ \text{21, Lemma 3.1}\right] $
we have

\textbf{Lemma A}$_{1}$. Assume $A$ is a $\varphi -$ positive linear operator
on a Banach space $E$. Then for any\ $h>0$ and $0\leq \mu \leq 1-\varkappa $
the operator-function \ 

\begin{equation*}
\Psi _{t}\left( \xi \right) =\Psi _{t,h,\mu }\left( \xi \right) =T\left(
t\right) \left\vert \xi \right\vert ^{\alpha +\sigma }A^{1-\varkappa -\mu
}h^{-\mu }\left[ A+\sum\limits_{k=1}^{n}t_{k}\left\vert \xi _{k}\right\vert
^{l_{k}}+h^{-1}\right] ^{-1}
\end{equation*}%
is bounded in\ $E$ uniformly with respect to $\xi \in R^{n}$,$\ h>0$ and $t,$
i.e. there exists a constant $C_{\mu }$ such that 
\begin{equation}
\ \ \ \left\Vert \Psi _{t,h,\mu }\left( \xi \right) \right\Vert _{L\left(
E\right) }\leq C_{\mu }  \tag{3.1}
\end{equation}%
for all $\xi \in R^{n}$ and $h>0.$

In a similar way as in $\left[ 19,\text{Theorem 2}\right] $ we obtain the
following

\textbf{Theorem A. }Assume$\ E_{0}$, $E$ are two Banach spaces and the
embedding $E_{0}\subset E$ is compact. Let $\Omega $ be a bounded domain in $%
R^{n}$ and $\gamma \in A_{p_{k}}$ for $p_{k}\in \left( 1,\infty \right) $.
Then the embedding 
\begin{equation*}
W_{\mathbf{p},\gamma }^{l}\left( \Omega ;E_{0},E\right) \subset L_{\mathbf{p}%
,\gamma }\left( \Omega ;E\right)
\end{equation*}%
is compact.

Let 
\begin{equation*}
X=L_{\mathbf{p},\gamma }\left( R^{n};E\right) ,\text{ }Y=W_{\mathbf{p,\gamma 
}}^{l}\left( R^{n};E\left( A\right) ,E\right) .
\end{equation*}

One of main result of this section is the following:

\ \ \ \ \ \ \textbf{Theorem 3.1}. Assume$\ E$ is a Banach space satisfying
the multiplier condition with respect to\ $\mathbf{p}$, $\mathbf{q}$ and
weighted function $\gamma $. Suppose $A$ is a $\varphi -$ positive operator
in $E$. Then for $0\leq \mu \leq 1-\varkappa $, $1<p_{k}\leq q_{k}<\infty $ $%
\ $or $0<\mu <1-\varkappa $ for $1\leq p_{k}\leq q_{k}\leq \infty $\ the
embedding 
\begin{equation*}
D^{\alpha }W_{\mathbf{p},\gamma }^{l}\left( R^{n};E\left( A\right) ,E\right)
\subset L_{\mathbf{q,\gamma }}\left( R^{n};E\left( A^{1-\varkappa -\mu
}\right) \right)
\end{equation*}%
is a continuous and there exists a constant $C_{\mu }$ \ $>0$ depending only
on $\mu $ such that 
\begin{equation}
T\left( t\right) \left\Vert D^{\alpha }u\right\Vert _{L_{\mathbf{q,\gamma }%
}\left( R^{n};E\left( A^{1-\varkappa -\mu }\right) \right) }\leq C_{\mu }%
\left[ h^{\mu }\left\Vert u\right\Vert _{W_{\mathbf{p},\gamma ,t}^{l}\left(
R^{n};E\left( A\right) ,E\right) }+h^{-\left( 1-\mu \right) }\left\Vert
u\right\Vert _{X}\right]  \tag{3.2}
\end{equation}%
for $u\in Y$ and $h>0.$

\ \ \textbf{Proof}. We have 
\begin{equation}
\left\Vert D^{\alpha }u\right\Vert _{L_{\mathbf{q},\gamma }\left(
R^{n};E\left( A^{1-\varkappa -\mu }\right) \right) }\backsim  \tag{3.3}
\end{equation}%
\begin{equation*}
\left( \left( \int\limits_{-\infty }^{\infty }\left( ...\int\limits_{-\infty
}^{\infty }\left( \int\limits_{-\infty }^{\infty }\left\Vert D^{\alpha
}u\right\Vert _{E\left( A^{1-\varkappa -\mu }\right) }^{q_{1}}\gamma \left(
x\right) dx_{1}\right) ^{\frac{q_{2}}{q_{1}}}dx_{2}\right) ^{\frac{q_{3}}{%
q_{2}}}...\right) ^{\frac{q_{n}}{q_{n-1}}}dx_{n}\right) ^{\frac{1}{q_{n}}%
}\sim
\end{equation*}%
\begin{equation*}
\left( \left( \int\limits_{-\infty }^{\infty }\left( ...\int\limits_{-\infty
}^{\infty }\left( \int\limits_{-\infty }^{\infty }\left\Vert A^{1-\varkappa
-\mu }D^{\alpha }u\right\Vert _{E}^{q_{1}}{}\gamma \left( x\right)
dx_{1}\right) ^{\frac{q_{2}}{q_{1}}}dx_{2}\right) ^{\frac{q_{3}}{q_{2}}%
}...\right) ^{\frac{q_{n}}{q_{n-1}}}dx_{n}\right) ^{\frac{1}{q_{n}}}
\end{equation*}%
for all\ $u$ such that 
\begin{equation*}
\left\Vert D^{\alpha }u\right\Vert _{L_{\mathbf{q},\gamma }\left(
R^{n};E\left( A^{1-\varkappa -\mu }\right) \right) }<\infty .
\end{equation*}%
On the other hand, it is clear to see that

\begin{equation}
A^{1-\alpha -\mu }D^{\alpha }u=F^{-1}FA^{1-\varkappa -\mu }D^{\alpha
}u=F^{-1}A^{1-\varkappa -\mu }FD^{\alpha }u=  \tag{3.4}
\end{equation}%
\begin{equation*}
F^{-1}A^{1-\varkappa -\mu }\left( i\xi \right) ^{\alpha }Fu=F^{-1}\left(
i\xi \right) ^{\alpha }A^{1-\varkappa -\mu }Fu.
\end{equation*}%
Hence, denoting \ $Fu$ by $\hat{u},$ we get from relations $\left(
3.3\right) $ and $\left( 3.4\right) $%
\begin{equation*}
\left\Vert D^{\alpha }u\right\Vert _{L_{\mathbf{q},\gamma }\left(
R^{n};E\left( A^{1-\varkappa -\mu }\right) \right) }\backsim \left\Vert
F^{-1}\left( i\xi \right) ^{\alpha }A^{1-\varkappa -\mu }\hat{u}\right\Vert
_{L_{\mathbf{q},\gamma }\left( R^{n};E\right) }.
\end{equation*}%
Similarly, in view of Definition 7 for $u\in Y$ we have 
\begin{equation*}
\left\Vert u\right\Vert _{W_{\mathbf{p},\gamma ,t}^{l},\left( R^{n};E\left(
A\right) ,E\right) }=\left\Vert u\right\Vert _{L_{\mathbf{p},\gamma }\left(
R^{n};E\left( A\right) \right) }+\sum\limits_{k=1}^{n}\left\Vert
t_{k}D_{k}^{l_{k}}u\right\Vert _{X}=
\end{equation*}%
\begin{equation*}
\left\Vert F^{-1}\hat{u}\right\Vert _{X}+\sum\limits_{k=1}^{n}\left\Vert
t_{k}F^{-1}\left[ \left( i\xi _{k}\right) ^{l_{k}}\hat{u}\right] \right\Vert
_{X}\backsim \left\Vert F^{-1}A\hat{u}\right\Vert
_{X}+\sum\limits_{k=1}^{n}\left\Vert t_{k}F^{-\shortmid }\left[ \left( i\xi
_{k}\right) ^{l_{k}}\hat{u}\right] \right\Vert _{X}.
\end{equation*}%
Therefore, for proving the inequality $\left( 3.2\right) $ it suffices to
show 
\begin{equation*}
T\left( t\right) \left\Vert F^{-1}\left( i\xi \right) ^{\alpha
}A^{1-\varkappa -\mu }\hat{u}\right\Vert _{L_{\mathbf{q},\gamma }\left(
R^{n};E\right) }\leq
\end{equation*}%
\begin{equation}
C_{\mu }(h^{\mu }\left\Vert F^{-1}A\hat{u}\right\Vert
_{X}+\sum\limits_{k=1}^{n}\left\Vert t_{k}F^{-1}\left[ \left( i\xi
_{k}\right) ^{l_{k}}\hat{u}\right] \right\Vert _{X}+h^{-\left( 1-\mu \right)
}\left\Vert F^{-1}\hat{u}\right\Vert _{X}).  \tag{3.5}
\end{equation}%
Therefore, the inequality $\left( 3.2\right) $ will follow if we prove the
following estimate \ 
\begin{equation}
T\left( t\right) \left\Vert F^{-1}\left[ \xi ^{\alpha +\sigma
}A^{1-\varkappa -\mu }\hat{u}\right] \right\Vert _{L_{\mathbf{q},\gamma
}\left( R^{n};E\right) }\leq C_{\mu }\left\Vert F^{-1}G\left( \xi \right) 
\hat{u}\right\Vert _{X}\text{.}  \tag{3.6}
\end{equation}%
for $u\in Y,$ where%
\begin{equation*}
G\left( \xi \right) =h^{\mu }\left[ A+\sum\limits_{k=1}^{n}t_{k}\left\vert
\xi _{k}\right\vert ^{l_{k}}+h^{-\left( 1-\mu \right) }\right] .
\end{equation*}%
\ Due to positivity of\ $A,$ the operator function $G\left( \xi \right) $
has a bounded inverse in\ $E$ for all $\xi \in R^{n}.$ So, we can set%
\begin{equation*}
T\left( t\right) F^{-1}\xi ^{\alpha +\sigma }A^{1-\varkappa -\mu }\hat{u}=
\end{equation*}%
\begin{equation}
T\left( t\right) F^{-1}\xi ^{\alpha +\sigma }A^{1-\varkappa -\mu
}G^{-1}\left( \xi \right) \left[ h^{\mu }\left(
A+\sum\limits_{k=1}^{n}t_{k}\left\vert \xi _{k}\right\vert ^{l_{k}}\right)
+h^{-\left( 1-\mu \right) }\right] \hat{u}.  \tag{3.7}
\end{equation}%
By Definition 6 it is clear to see that the inequality $\left( 3.6\right) $
will follow immediately from $\left( 3.7\right) $ if we can prove that the
operator-function

\begin{equation*}
\Psi _{t,h,\mu }=\xi ^{\alpha +\sigma }A^{1-\varkappa -\mu }\left[ h^{\mu
}(A+\sum\limits_{k=1}^{n}t_{k}\left\vert \xi _{k}\right\vert
^{l_{k}})+h^{-\left( 1-\mu \right) }\right] ^{-1}
\end{equation*}%
is a multiplier in $M_{\mathbf{p},\gamma }^{\mathbf{q},\gamma }\left(
E\right) $ uniformly with respect to\ $h$ and $t.$ So, it suffices to show
that for all $\beta \in U_{n}$ and $\xi \in V_{n}$\ there exists a constant
\ $C_{\mu }>0$ such that the following uniform estimate holds%
\begin{equation}
\left\Vert D^{\beta }\Psi _{t,h,\mu }\left( \xi \right) \right\Vert
_{L\left( E\right) }\leq C_{\mu }\dprod\limits_{k=1}^{n}\left\vert \xi
_{k}\right\vert ^{-\left( \beta _{k}+\sigma _{k}\right) }.  \tag{3.8}
\end{equation}%
\ To see this, by applying the Lemma A$_{1}$ for all $\xi \in R^{n}$ we get
a constant $C_{\mu }>0$ depending only on $\mu $ such that 
\begin{equation}
\ \left\Vert \Psi _{t,h,\mu }\left( \xi \right) \right\Vert _{L\left(
E\right) }\leq C_{\mu }\ \dprod\limits_{k=1}^{n}\left\vert \xi
_{k}\right\vert ^{-\sigma _{k}}.\ \ \ \ \ \ \ \ \   \tag{3.9}
\end{equation}%
This shows that the inequality $\left( 3.8\right) $ is satisfies for $\beta
=\left( 0,...,0\right) .$ Now, we next consider $\left( 3.8\right) $ for $%
\beta =\left( \beta _{1},...\beta _{n}\right) $ where \ $\beta _{k}=1$ and $%
\beta =0$ for $j\neq k$. Then, by using the positivity properties of $A$ we
obtain

\begin{equation*}
\left\vert \frac{\partial }{\partial \xi _{k}}\Psi _{t}\left( \xi \right)
\right\vert \leq \prod\limits_{k=1}^{n}t_{k}^{\frac{\alpha _{k}}{l_{k}}%
}\left( i\right) ^{\left\vert \alpha \right\vert }\alpha _{k}\left\vert \xi
_{1}^{\alpha _{1}}...\xi _{k-1}^{\alpha _{k-1}}\xi _{k}^{\alpha
_{k}-1}...\xi _{n}^{\alpha _{n}}\right\vert
\end{equation*}%
\begin{equation*}
\left\Vert A^{1-\varkappa -\mu }\left[ h^{\mu }\left(
A+\sum\limits_{k=1}^{n}t_{k}\left\vert \xi _{k}\right\vert ^{l_{k}}\right)
+h^{-\left( 1-\mu \right) }\right] ^{-1}\right\Vert +
\end{equation*}%
\begin{equation*}
\left\vert \xi \right\vert ^{\alpha }\left\Vert A^{1-\varkappa -\mu }\left[
h^{\mu }\left( A+\sum\limits_{k=1}^{n}t_{k}\left\vert \xi _{k}\right\vert
^{l_{k}}\right) +h^{-\left( 1-\mu \right) }\right] ^{-2}\right\Vert
ht_{k}\left\vert \xi _{k}\right\vert ^{l_{k}-1}\leq
\end{equation*}%
\begin{equation*}
C_{\mu }\left\vert \xi _{k}\right\vert
^{-1}\dprod\limits_{j=1}^{n}\left\vert \xi _{j}\right\vert ^{-\sigma _{j}},%
\text{ }k=1,2...n.
\end{equation*}%
Repeating the above process, we obtain that for all $\beta \in U_{n},$ $\xi
\in V_{n}$ there exists a constant \ $C_{\mu }>0$ \ depending only $\mu $
such that 
\begin{equation*}
\ \ \ \ \ \ \ \ \ \ \ \left\Vert D^{\beta }\Psi _{t}\left( \xi \right)
\right\Vert _{L\left( E\right) }\leq C_{\mu }\left\vert \xi \right\vert
^{-\left( \beta +\sigma \right) }.
\end{equation*}%
\ Therefore, the operator-function $\Psi _{t,h,\mu }\left( \xi \right) $ is
a uniform multiplier with respect to $h$ and $t,$ i.e, 
\begin{equation*}
\Psi _{t,h,\mu }\in M_{\mathbf{p},\gamma }^{\mathbf{q},\gamma }\left(
E\right) ,\text{ for }t_{k},\text{ }h,\text{ }\mu \in \text{\ }\mathbb{R}%
_{+}.
\end{equation*}

This completes the proof of the Theorem 3.1.

\ \ It is possible to state Theorem 3.1 in a more general setting. For this
aim, we use the concept of extension operator.

\ \textbf{Condition 3.1}. Let $A$ be positive operator in Banach spaces $E$
satisfying multiplier condition with respect to\textbf{\ }$\mathbf{p}$ and
weighted function $\gamma .$\ Assume a region \ $\Omega \subset R^{n}$ such
that there exists bounded linear extension operator$\ B$ from $W_{\mathbf{p}%
,\gamma }^{l}\left( \Omega ,E\left( A\right) ,E\right) $ to $W_{\mathbf{p}%
,\gamma }^{l}\left( R^{n},E\left( A\right) ,E\right) $ for \ $1\leq
p_{k}\leq \infty .$

\textbf{Remark 3.1}. If \ $\Omega \subset R^{n}$ is a region satisfying the
strong \ $l-$horn condition \ (see $\left[ 3\right] $, p.117), $E=\mathbb{C}%
, $ $A$ $=I$ \ and $\gamma \left( x\right) \equiv 1$\ then for $1<p<\infty $
there exists a bounded linear extension operator from \ $W_{p}^{l}\left(
\Omega \right) =W_{p}^{l}\left( \Omega ;\mathbb{R},\mathbb{R}\right) $ to $%
W_{p}^{l}\left( R^{n}\right) =W_{p}^{l}\left( R^{n};\mathbb{R},\mathbb{R}%
\right) .$

\textbf{Theorem 3.2}. Assume conditions of Theorem 3.1 and Condition 3.1 are
hold. Then for $0\leq \mu \leq 1-\varkappa $ the embedding 
\begin{equation*}
D^{\alpha }W_{\mathbf{p},\gamma }^{l}\left( \Omega ;E\left( A\right)
,E\right) \subset L_{\mathbf{q},\gamma }\left( \Omega ;E\left(
A^{1-\varkappa -\mu }\right) \right)
\end{equation*}%
is continuous and there exists a constant $C_{\mu }$ depending only on $\mu $
such that 
\begin{equation}
T\left( t\right) \left\Vert D^{\alpha }u\right\Vert _{L_{\mathbf{q},\gamma
}\left( \Omega ;E\left( A^{1-\varkappa -\mu }\right) \right) }\leq 
\tag{3.10}
\end{equation}%
\begin{equation*}
C_{\mu }\left[ h^{\mu }\left\Vert u\right\Vert _{W_{\mathbf{p},\gamma
,t}^{l}\left( \Omega ;E\left( A\right) ,E\right) }+h^{-\left( 1-\mu \right)
}\left\Vert u\right\Vert _{L_{\mathbf{p},\gamma }\left( \Omega ;E\right) }%
\right]
\end{equation*}%
for $u\in W_{\mathbf{p},\gamma }^{l}\left( \Omega ;E\left( A\right)
,E\right) $ and \ $h>0.$

\ \textbf{Proof}.\ It is suffices to prove the estimate $\left( 3.10\right)
. $ Let $B$ is a bounded linear extension operator\ from \ $W_{\mathbf{p}%
,\gamma }^{l}\left( \Omega ;E\left( A\right) ,E\right) $\ to \ $W_{\mathbf{p}%
,\gamma }^{l}\left( R^{n};E\left( A\right) ,E\right) ,$ and let $B_{\Omega }$%
\ be\ the restriction operator from \ $R^{n}$ to $\Omega .$ Then for any $%
u\in W_{\mathbf{p},\gamma }^{l}\left( \Omega ;E\left( A\right) ,E\right) $\
we have 
\begin{equation*}
T\left( t\right) \left\Vert D^{\alpha }u\right\Vert _{L_{\mathbf{q},\gamma
}\left( \Omega ;E\left( A^{1-\varkappa -\mu }\right) \right) }=T\left(
t\right) \left\Vert D^{\alpha }B_{\Omega }Bu\right\Vert _{L_{\mathbf{q}%
,\gamma }\left( \Omega ;E\left( A^{1-\varkappa -\mu }\right) \right) }\leq
\end{equation*}%
\begin{equation*}
CT\left( t\right) \left\Vert D^{\alpha }Bu\right\Vert _{L_{\mathbf{q},\gamma
}\left( R^{n};E\left( A^{1-\varkappa -\mu }\right) \right) }\leq
\end{equation*}%
\begin{equation*}
C_{\mu }\left[ h^{\mu }\left\Vert Bu\right\Vert _{W_{\mathbf{p},\gamma
,t}^{l}\left( R^{n};E\left( A\right) ,E\right) }+h^{-\left( 1-\mu \right)
}\left\Vert Bu\right\Vert _{L_{\mathbf{p},\gamma }\left( R^{n};E\right) }%
\right] \leq
\end{equation*}%
\begin{equation*}
C_{\mu }\left[ h^{\mu }\left\Vert u\right\Vert _{W_{\mathbf{p},\gamma
,t}^{l}\left( \Omega ;E\left( A\right) E\right) }+h^{-\left( 1-\mu \right)
}\left\Vert u\right\Vert _{L_{\mathbf{p},\gamma }\left( \Omega ;E\right) }%
\right] .
\end{equation*}%
\textbf{Result 3.1}. Assume the conditions of Theorem 3.2 are satisfied.
Then for $u\in W_{\mathbf{p},\gamma }^{l}\left( \Omega ;E\left( A\right)
,E\right) $ we have the following multiplicative estimate 
\begin{equation}
\left\Vert D^{\alpha }u\right\Vert _{L_{\mathbf{q},\gamma }\left( \Omega
;E\left( A^{1-\varkappa -\mu }\right) \right) }\leq C_{\mu }\left\Vert
u\right\Vert _{W_{\mathbf{p},\gamma }^{l}\left( \Omega ;E\left( A\right)
E\right) }^{1-\mu }.\left\Vert u\right\Vert _{L_{\mathbf{p},\gamma }\left(
\Omega ;E\right) }^{\mu }.  \tag{3.11}
\end{equation}%
\ Indeed, setting 
\begin{equation*}
h=\left\Vert u\right\Vert _{L_{\mathbf{p},\gamma }\left( \Omega ;E\right)
}.\left\Vert u\right\Vert _{W_{\mathbf{p},\gamma }^{l}\left( \Omega ;E\left(
A\right) ,E\right) }^{-1}
\end{equation*}%
\ in\ $\left( 3.10\right) $ we obtain $\left( 3.11\right) .$

\textbf{Theorem 3.3. }Assume that the conditions of Theorem 3.2 are
satisfied. Suppose $\Omega $ is a bounded domain in $R^{n}$ and $A^{-1}$ is
a compact operator in $E.$ Let $\gamma \in A_{p_{k}}$ for $p_{k}\in \left(
1,\infty \right) $. Then for $0<\mu \leq 1-\varkappa $ the embedding%
\begin{equation*}
D^{\alpha }W_{\mathbf{p},\gamma }^{l}\left( \Omega ;E\left( A\right)
,E\right) \subset L_{\mathbf{q},\gamma }\left( \Omega ;E\left(
A^{1-\varkappa -\mu }\right) \right)
\end{equation*}%
is compact.

\textbf{Proof. }By virtue of Theorem $A$ we get that the embedding%
\begin{equation*}
W_{\mathbf{p},\gamma }^{l}\left( \Omega ;E\left( A\right) ,E\right) \subset
L_{\mathbf{p},\gamma }\left( \Omega ;E\right)
\end{equation*}%
is compact. Then by $\left( 3.11\right) $ we obtain the assertion of Theorem
3.3.

Let 
\begin{equation*}
\sigma =\max\limits_{k}p_{k}.
\end{equation*}

\bigskip \textbf{Theorem 3}.\textbf{4.} Suppose conditions of Theorem 3.1
are hold$.$ Then for $0<\mu <1-\varkappa $ the embedding 
\begin{equation*}
D^{\alpha }W_{\mathbf{p},\gamma }^{l}\left( R^{n};E\left( A\right) ,E\right)
\subset L_{\mathbf{q},\gamma }\left( R^{n};\left( E\left( A\right) ,E\right)
_{\varkappa +\mu ,\sigma }\right)
\end{equation*}%
is continuous and there exists a constant $C_{\mu }$ depending only on $\mu $
such that

\begin{equation}
\left\Vert D^{\alpha }u\right\Vert _{L_{\mathbf{q},\gamma }\left(
R^{n};\left( E\left( A\right) ,E\right) _{\varkappa +\mu ,\sigma }\right)
}\leq h^{\mu }\left\Vert u\right\Vert _{W_{\mathbf{p},\gamma }^{l}\left(
R^{n};E\left( A\right) ,E\right) }+h^{-\left( 1-\mu \right) }\left\Vert
u\right\Vert _{L_{\mathbf{p},\gamma }\left( R^{n};E\right) }  \tag{3.12}
\end{equation}%
for $u\in Y$ and $0<h\leq h_{0}<\infty .$

\textbf{Proof. }It is sufficient to prove the estimate $\left( 3.12\right) $
for $u\in Y.$ By definition of interpolation spaces $\left( E\left( A\right)
,E\right) _{\varkappa +\mu ,\sigma }$ (see $\left[ \text{22, \S 1.14.5}%
\right] $) the estimate $\left( 3.12\right) $ is equivalent to the inequality

\begin{equation}
\left\Vert F^{-1}y^{1-\varkappa -\mu -\frac{1}{\sigma }}\left[ A^{\chi +\mu
}\left( A+y\right) ^{-1}\right] \xi ^{\alpha }\hat{u}\right\Vert _{L_{%
\mathbf{p},\gamma }\left( R_{+}^{n+1};E\right) }  \tag{3.13}
\end{equation}

\begin{equation*}
\leq C_{\mu }\left\Vert F^{-1}\left[ h^{\mu
}(A+\sum\limits_{k=1}^{n}A+\sum\limits_{k=1}^{n}\left\vert \xi
_{k}\right\vert ^{l_{k}}+h^{-\left( 1-\mu \right) }\right] \hat{u}%
\right\Vert _{L_{\mathbf{p},\gamma }\left( R^{n};E\right) }.
\end{equation*}%
By multiplier properties, the inequality $\left( 3.13\right) $ will follow
immediately if we will prove that the operator-function

\begin{equation*}
\Psi =\left( i\xi \right) ^{\alpha }y^{1-\varkappa -\mu -\frac{1}{p}}A^{\chi
+\mu }\left( A+y\right) ^{-1}\left[ h^{\mu }\left(
A+\sum\limits_{k=1}^{n}\left\vert \xi _{k}\right\vert ^{l_{k}}\right)
+h^{-\left( 1-\mu \right) }\right] ^{-1}
\end{equation*}%
is a multiplier from $L_{\mathbf{p},\gamma }\left( R^{n};E\right) $ to $L_{%
\mathbf{p},\gamma }\left( R^{n};L_{\sigma }\left( R_{+};E\right) \right) .$
This fact is proved by the same manner as Theorem 3.1. Therefore, we get the
estimate $\left( 3.12\right) .$

In a similar way, as the Theorem 3.2 we obtain

\textbf{Theorem 3.5}. Suppose conditions of Theorem 3.2 are hold$.$ Then for 
$0<\mu <1-\varkappa $ the embedding 
\begin{equation*}
D^{\alpha }W_{\mathbf{p},\gamma }^{l}\left( \Omega ;E\left( A\right)
,E\right) \subset L_{\mathbf{q},\gamma }\left( \Omega ;\left( E\left(
A\right) ,E\right) _{\varkappa +\mu ,\sigma }\right)
\end{equation*}%
is continuous and there exists a constant $C_{\mu }$ depending only on $\mu $
such that 
\begin{equation}
\left\Vert D^{\alpha }u\right\Vert _{L_{\mathbf{q},\gamma }\left( \Omega
,\left( E\left( A\right) ,E\right) _{\varkappa +\mu ,p}\right) }\leq C_{\mu }%
\left[ h^{\mu }\left\Vert u\right\Vert _{W_{\mathbf{p},\gamma }^{l}\left(
\Omega ;E\left( A\right) ,E\right) }+h^{-\left( 1-\mu \right) }\left\Vert
u\right\Vert _{L_{\mathbf{p},\gamma }\left( \Omega ;E\right) }\right] 
\tag{3.14}
\end{equation}%
for $u\in W_{\mathbf{p},\gamma }^{l}\left( \Omega ;E\left( A\right)
,E\right) $ and $0<h\leq h_{0}<\infty .$

\textbf{Result 3. 2}. Suppose the conditions of Theorem 3.2 are hold. Then
for $u\in W_{\mathbf{p},\gamma }^{l}\left( \Omega ;E\left( A\right)
,E\right) $\ we have the following multiplicative estimate 
\begin{equation}
\left\Vert D^{\alpha }u\right\Vert _{L_{\mathbf{q},\gamma }\left( \Omega
;\left( E\left( A\right) ,E\right) _{\varkappa +\mu ,p}\right) }\leq C_{\mu
}\left\Vert u\right\Vert _{W_{\mathbf{p},\gamma }^{l}\left( \Omega ;E\left(
A\right) ,E\right) }^{1-\mu }\left\Vert u\right\Vert _{L_{\mathbf{p},\gamma
}\left( \Omega ;E\right) }^{\mu }.  \tag{3.15}
\end{equation}%
Indeed setting $\left\Vert u\right\Vert _{L_{\mathbf{p},\gamma }\left(
\Omega ;E\right) }.\left\Vert u\right\Vert _{W_{\mathbf{p},\gamma
}^{l}\left( \Omega ;E\left( A\right) ,E\right) }^{-1}$ in $\left(
3.14\right) $ we obtain\ $\left( 3.15\right) .$

From the estimate $\left( 3.15\right) $ and Theorem A, in a similar way as
Theorem 3.3 we obtain

\textbf{Theorem 3.6. }Assume that the conditions of Theorem 3.2 are
satisfied. Suppose $\Omega $ is a bounded domain in $R^{n}$ and $A^{-1}$ is
a compact operator in $E.$ Let $\gamma \in A_{p_{k}}$ for $p_{k}\in \left(
1,\infty \right) $. Then for $0<\mu \leq 1-\varkappa $ the embedding%
\begin{equation*}
D^{\alpha }W_{\mathbf{p},\gamma }^{l}\left( \Omega ;E\left( A\right)
,E\right) \subset L_{\mathbf{q},\gamma }\left( \Omega ;\left( E\left(
A\right) ,E\right) _{\varkappa +\mu ,\sigma }\right)
\end{equation*}

is compact.

From Theorem 3.2 we obtain

\textbf{Result 3.2}. Assume the conditions of Theorem 3.2 are satisfied for $%
l_{1}=l_{2}=\ldots =l_{n}=m.$ Then for $0\leq \mu \leq 1-\varkappa $ the
embedding 
\begin{equation*}
D^{\alpha }W_{\mathbf{p},\gamma }^{m}\left( \Omega ;E\left( A\right)
,E\right) \subset L_{\mathbf{q},\gamma }\left( \Omega ;E\left(
A^{1-\varkappa -\mu }\right) \right)
\end{equation*}%
is continuous and there exists a constant $C_{\mu }$ depending only on $\mu $
such that 
\begin{equation*}
T\left( t\right) \left\Vert D^{\alpha }u\right\Vert _{L_{\mathbf{q},\gamma
}\left( \Omega ;E\left( A^{1-\varkappa -\mu }\right) \right) }\leq C_{\mu }%
\left[ h^{\mu }\left\Vert u\right\Vert _{W_{\mathbf{p},\gamma ,t}^{m}\left(
\Omega ;E\left( A\right) ,E\right) }+h^{-\left( 1-\mu \right) }\left\Vert
u\right\Vert _{L_{\mathbf{p},\gamma }\left( \Omega ;E\right) }\right]
\end{equation*}%
for $u\in W_{\mathbf{p},\gamma }^{m}\left( \Omega ;E\left( A\right)
,E\right) $ and\ $h>0,$ where 
\begin{equation*}
\varkappa =\frac{1}{m}\left( \left\vert \alpha \right\vert
+\dsum\limits_{k=1}^{n}\sigma _{k}\right) .
\end{equation*}

\textbf{Result 3.3}. If $E=H$, where $H$ is a Hilbert space and \ $%
p_{k}=q_{k}=2,$ $\Omega =\left( 0,T\right) ,$ $l_{1}=l_{2}=\ldots =l_{n}=m$ $%
,$ $A=A^{\times }\geq cI,$\ $\gamma \left( x\right) \equiv 1$\ then we
obtain the well known Lions-Peetre $\left[ 11\right] $ result. \ Moreover,
the result of Lions-Peetre are improving even in the one dimensional case
for the non selfedjoint positive operators $A.$

From Theorems 3.2, 3.3 we obtain

\textbf{Result 3.4}. Suppose the conditions of Theorem 3.2 are satisfied for 
$\gamma \left( x\right) \equiv 1.$ Then for $0\leq \mu \leq 1-\varkappa $
the embedding 
\begin{equation*}
D^{\alpha }W_{\mathbf{p}}^{l}\left( \Omega ;E\left( A\right) ,E\right)
\subset L_{\mathbf{q}}\left( \Omega ;E\left( A^{1-\varkappa -\mu }\right)
\right)
\end{equation*}%
is a continuous and there exists a constant $C_{\mu }$ \ $>0$, depending
only on $\mu $ such that 
\begin{equation*}
T\left( t\right) \left\Vert D^{\alpha }u\right\Vert _{L_{\mathbf{q}}\left(
\Omega ;E\left( A^{1-\varkappa -\mu }\right) \right) }\leq C_{\mu }\left[
h^{\mu }\left\Vert u\right\Vert _{W_{\mathbf{p},t}^{l}\left( \Omega ;E\left(
A\right) ,E\right) }+h^{-\left( 1-\mu \right) }\left\Vert u\right\Vert _{L_{%
\mathbf{p}}\left( \Omega ;E\right) }\right]
\end{equation*}%
for $u\in W_{\mathbf{p}}^{l}\left( \Omega ;E\left( A\right) ,E\right) $ and $%
h>0$.

Moreover, if $\Omega $ is a bounded domain in $R^{n}$ and $A^{-1}$ is a
compact operator in $E,$ then for $0<\mu \leq 1-\varkappa $ the embedding%
\begin{equation*}
D^{\alpha }W_{\mathbf{p}}^{l}\left( \Omega ;E\left( A\right) ,E\right)
\subset L_{\mathbf{q}}\left( \Omega ;E\left( A^{1-\varkappa -\mu }\right)
\right)
\end{equation*}%
is compact.

If \ $E=\mathbb{C},$ $A=I$, $\gamma \left( x\right) \equiv 1$ we get the
embedding $D^{\alpha }W_{\mathbf{p}}^{l}\left( \Omega \right) \subset L_{%
\mathbf{q}}\left( \Omega \right) $ proved in $\left[ 3\right] $ for
numerical Sobolev spaces $W_{\mathbf{p}}^{l}\left( \Omega \right) .$

\ \ \ \ \ \ \ \ \ \ \ \ \ \ \ \ \ \ \ \ \ \ \ \ \ \ \ \ \ \ \ \ \ \ \ 

\begin{center}
\ \ \ \ \ \ \ \ \ \ \ \ \ \ \ \ \ \ \ \ \ \ \ \ \ \ \ \ \ \ \ \ \ \ \ 
\textbf{4. Application}
\end{center}

\ Let\ $s>0.$ Consider the following sequence space (see e.g. $\left[ \text{%
22, \S\ 1.18}\right] $)%
\begin{equation*}
l_{q}^{s}=\left\{ u=\left\{ u_{i}\right\} ,\text{ }i=1,2,...,\infty ,\text{ }%
u_{i}\in \mathbb{C}\right\}
\end{equation*}%
\ with the norm$\ $%
\begin{equation*}
\ \ \ \ \ \ \left\Vert u\right\Vert _{l_{\nu }^{s}}=\left(
\sum\limits_{i=1}^{\infty }2^{i\nu s}\left\vert u_{i}\right\vert ^{p}\right)
^{\frac{1}{q}}<\infty ,\text{ }\nu \in \left( 1,\infty \right) .
\end{equation*}%
Note that, \ $l_{\nu }^{0}=l_{\nu }.$ Let $A$ be infinite matrix defined in $%
l_{\nu }$ such that \ $D\left( A\right) =l_{\nu }^{s},$ $A=\left[ \delta
_{ij}2^{si}\right] ,\ $where \ $\delta _{ij}=0$, when \ $i\neq j,$ \ $\delta
_{ij}=1,$ when $i=j=1,2,...,\infty .$

It is clear to see that, the operator $A$\ is positive in\ $l_{\nu }.$Then
from Theorem 3.2 and Theorem 3.3 we obtain the following results

\textbf{Result 4.1. }Suppose the conditions of Theorem 3.2 are satisfied for 
$E=\mathbb{C}$. Then for $0\leq \mu \leq 1-\varkappa $, $1<p_{k}\leq
q_{k}<\infty $ $\ $or $0<\mu <1-\varkappa $ for $1\leq p_{k}\leq q_{k}\leq
\infty $ the embedding 
\begin{equation*}
D^{\alpha }W_{\mathbf{p},\gamma }^{l}\left( \Omega ,l_{\nu }^{s},l_{\nu
}\right) \subset L_{\mathbf{q},\gamma }\left( \Omega ,l_{\nu }^{s\left(
1-\varkappa -\mu \right) }\right)
\end{equation*}%
is a continuous and there exists a constant $C_{\mu }$ \ $>0$, depending
only on $\mu $ such that 
\begin{equation*}
T\left( t\right) \left\Vert D^{\alpha }u\right\Vert _{L_{\mathbf{q,\gamma }%
}\left( \Omega ;l_{\nu }^{s\left( 1-\varkappa -\mu \right) }\right) }\leq
C_{\mu }\left[ h^{\mu }\left\Vert u\right\Vert _{W_{\mathbf{p},\gamma
,t}^{l}\left( \Omega ;l_{\nu }^{s},l_{\nu }\right) }+h^{-\left( 1-\mu
\right) }\left\Vert u\right\Vert _{L_{\mathbf{p,\gamma }}\left( \Omega
;l_{\nu }\right) }\right]
\end{equation*}%
for\ $u\in W_{\mathbf{p},\gamma }^{l}\left( \Omega ,l_{\nu }^{s},l_{\nu
}\right) $ and $h>0.$

\textbf{Result 4.2. }Suppose the conditions of Theorem 3.3 are hold for $E=%
\mathbb{C}$. Then for $0<\mu \leq 1-\varkappa ,$ $1<p_{k}\leq q_{k}<\infty $ 
$\ $or $0<\mu <1-\varkappa $ for $1\leq p_{k}\leq q_{k}\leq \infty $ the
embedding 
\begin{equation*}
D^{\alpha }W_{\mathbf{p},\gamma }^{l}\left( \Omega ,l_{\nu }^{s},l_{\nu
}\right) \subset L_{\mathbf{q},\gamma }\left( \Omega ,l_{\nu }^{s\left(
1-\varkappa -\mu \right) }\right)
\end{equation*}%
is compact.

\textbf{Result 4.3. }For $0\leq \mu \leq 1-\varkappa $, $1<p_{k}\leq
q_{k}<\infty $ $\ $or $0<\mu <1-\varkappa $ for $1\leq p_{k}\leq q_{k}\leq
\infty $ the embedding 
\begin{equation*}
D^{\alpha }W_{\mathbf{p}}^{l}\left( \Omega ,l_{\nu }^{s},l_{\nu }\right)
\subset L_{\mathbf{q}}\left( \Omega ,l_{\nu }^{s\left( 1-\varkappa -\mu
\right) }\right)
\end{equation*}%
is a continuous and there exists a constant $C_{\mu }$ \ $>0$, depending
only on $\mu $ such that 
\begin{equation*}
T\left( t\right) \left\Vert D^{\alpha }u\right\Vert _{L_{\mathbf{q}}\left(
\Omega ;l_{\nu }^{s\left( 1-\varkappa -\mu \right) }\right) }\leq C_{\mu }%
\left[ h^{\mu }\left\Vert u\right\Vert _{W_{\mathbf{p},t}^{l}\left( \Omega
;l_{\nu }^{s},l_{\nu }\right) }+h^{-\left( 1-\mu \right) }\left\Vert
u\right\Vert _{L_{\mathbf{p}}\left( \Omega ;l_{\nu }\right) }\right]
\end{equation*}%
for\ $u\in W_{\mathbf{p}}^{l}\left( \Omega ,l_{\nu }^{s},l_{\nu }\right) $
and $h>0.$

\textbf{Result 4.4. }For $0<\mu \leq 1-\varkappa ,$ $1<p_{k}\leq
q_{k}<\infty $ $\ $or $0<\mu <1-\varkappa $ for $1\leq p_{k}\leq q_{k}\leq
\infty $ the embedding 
\begin{equation*}
D^{\alpha }W_{\mathbf{p}}^{l}\left( \Omega ,l_{\nu }^{s},l_{\nu }\right)
\subset L_{\mathbf{q}}\left( \Omega ,l_{\nu }^{s\left( 1-\varkappa -\mu
\right) }\right)
\end{equation*}%
is compact.

Note that, these results haven't been obtained with classical method until
now.

\bigskip

\begin{center}
\textbf{\ \ \ \ \ \ \ \ \ \ \ \ \ \ \ \ \ \ \ \ \ \ \ \ 5. Separable
degenerate abstract differential operators }
\end{center}

\bigskip\ Let us consider the problem%
\begin{equation}
\sum\limits_{k=1}^{n}(-1)^{l_{k}}t_{k}D_{k}^{\left[ 2l_{k}\right] }u\left(
x\right) +\left( A+\lambda \right) u\left( x\right) +  \tag{5.1}
\end{equation}%
\begin{equation*}
\sum\limits_{\left\vert \alpha :2l\right\vert
<1}\prod\limits_{k=1}^{n}t_{k}^{\frac{\alpha _{k}}{2l_{k}}}A_{\alpha }(x)D^{%
\left[ \alpha \right] }u\left( x\right) =f\left( x\right) ,\text{ }x\in
R^{n},
\end{equation*}%
considered in $L_{\mathbf{p}}\left( R^{n};E\right) $, where\ $A$, $A_{\alpha
}$\ are linear operators in a Banach space $E$, $t_{k}$ are positive and $%
\lambda $ is a complex parameter.

Let 
\begin{equation*}
X=L_{\mathbf{p},\tilde{\gamma}}\left( R^{n};E\right) ,\text{ }Y=W_{\mathbf{p,%
\tilde{\gamma}}}^{2l}\left( R^{n};E\left( A\right) ,E\right) .
\end{equation*}

\textbf{Remark 5.1. }\bigskip Under the substitution 
\begin{equation}
\tau _{k}=\int\limits_{0}^{x_{k}}\gamma _{k}^{-1}\left( y\right) dy 
\tag{5.2}
\end{equation}%
the spaces\ $L_{\mathbf{p}}\left( R^{n};E\right) $, $W_{\mathbf{p},\gamma }^{%
\left[ 2l\right] }\left( R^{n};E\left( A\right) ,E\right) $\ are mapped
isomorphically onto the weighted spaces $X$ and $Y,$ where%
\begin{equation*}
\tilde{\gamma}=\tilde{\gamma}\left( \tau \right)
=\prod\limits_{k=1}^{n}\gamma _{k}\left( x_{k}\left( \tau _{k}\right)
\right) \text{, }\tau =\left( \tau _{1},\tau _{2},...,\tau _{n}\right) .
\end{equation*}

Moreover, under this transformation the problem $\left( 5.1\right) $ is
mapped to the following undegenerate problem 
\begin{equation*}
\sum\limits_{k=1}^{n}(-1)^{l_{k}}t_{k}D_{k}^{2l_{k}}\tilde{u}\left( \tau
\right) +\left( A+\lambda \right) \tilde{u}\left( \tau \right) +
\end{equation*}%
\begin{equation*}
\sum\limits_{\left\vert \alpha :2l\right\vert
<1}\prod\limits_{k=1}^{n}t_{k}^{\frac{\alpha _{k}}{2l_{k}}}A_{\alpha
}(x)D^{\alpha }\tilde{u}\left( \tau \right) =\tilde{f}\left( \tau \right)
,\tau \in R^{n},
\end{equation*}%
considered in the weighted space $L_{\mathbf{p},\tilde{\gamma}}\left(
R^{n};E\right) $, where 
\begin{equation*}
\tilde{u}\left( \tau \right) =u\left( x_{1}\left( \tau _{1}\right)
,x_{2}\left( \tau _{2}\right) ,...,x_{n}\left( \tau _{n}\right) \right) ,%
\text{ }\tilde{f}\left( \tau \right) =f\left( x_{1}\left( \tau _{1}\right)
,x_{2}\left( \tau _{2}\right) ,...,x_{n}\left( \tau _{n}\right) \right) .
\end{equation*}%
By redenoting $u=\tilde{u}\left( \tau \right) $ and $f=\tilde{f}\left( \tau
\right) $ we get 
\begin{equation}
\sum\limits_{k=1}^{n}(-1)^{l_{k}}t_{k}D_{k}^{2l_{k}}u\left( \tau \right)
+\left( A+\lambda \right) u\left( \tau \right) +  \tag{5.3}
\end{equation}

\begin{equation*}
\sum\limits_{\left\vert \alpha :2l\right\vert
<1}\prod\limits_{k=1}^{n}t_{k}^{\frac{\alpha _{k}}{2l_{k}}}A_{\alpha
}(x)D^{\alpha }u\left( \tau \right) =f\left( \tau \right) ,\text{ }\tau \in
R^{n}.
\end{equation*}

Consider first of all, the pr\i ncipal part of $\left( 5.3\right) $, i.e.
consider the problem 
\begin{equation}
\sum\limits_{k=1}^{n}(-1)^{l_{k}}t_{k}D_{k}^{2l_{k}}u\left( \tau \right)
+\left( A+\lambda \right) u\left( \tau \right) =f\left( \tau \right) ,\text{ 
}\tau \in R^{n}.  \tag{5.4}
\end{equation}

\textbf{Theorem 5.1.} Assume the following conditions are satisfied:

(1) $t_{k}>0,$ $\mathbf{p}=\left( p_{1},p_{2},...,p_{n}\right) ,$ $%
1<p_{k}<\infty $, $k=1,2,...,n$;

(2) $E$ is Banach space satisfying multiplier condition with respect to $%
\mathbf{p}$ and weighted function $\tilde{\gamma}$;

(3) $A$ is a $\varphi $-positive operator in Banach space $E$ for $0\leq
\varphi <\pi .$

Then for $\ f\in X$ and $\lambda \in S\left( \varphi \right) $ problem $%
\left( 5.4\right) $ has a unique solution $u\left( x\right) $ that belongs
to $Y$\ and the uniform coercive estimate holds 
\begin{equation}
\sum\limits_{k=1}^{n}\sum\limits_{i=0}^{2l_{k}}t_{k}^{\frac{i}{2l_{k}}%
}\left\vert \lambda \right\vert ^{1-\frac{i}{2l_{k}}}\left\Vert D_{k}^{\left[
i\right] }u\right\Vert _{X}+\left\Vert Au\right\Vert _{X}\leq C\left\Vert
f\right\Vert _{X}.  \tag{5.5}
\end{equation}%
\ \textbf{Proof.} \ By appllying Fourier transform to the equation $\left(
5.4\right) $ we\ obtain 
\begin{equation}
\sum\limits_{k=1}^{n}t_{k}\xi _{k}^{2l_{k}}\hat{u}\left( \xi \right) +\left(
A+\lambda \right) \hat{u}\left( \xi \right) =\hat{f}\left( \xi \right) . 
\tag{5.6}
\end{equation}%
\ 

It is clear that \ $\sum\limits_{k=1}^{n}t_{k}\xi _{k}^{2l_{k}}\geq 0$ \ for
all \ $\xi =\left( \xi _{1},...,\xi _{n}\right) \in R^{n}$. Therefore, we
get that \ $\lambda +\sum\limits_{k=1}^{n}t_{k}\xi _{k}^{2l_{k}}\in S\left(
\varphi \right) $ for all $\xi \in R^{n}$. Since $A$ is $\varphi $-positive,
we deduce that the operator function%
\begin{equation*}
\Psi \left( \xi \right) =A+\lambda +\sum\limits_{k=1}^{n}t_{k}\xi
_{k}^{2l_{k}}
\end{equation*}%
has a bounded inverse $\Psi ^{-1}\left( \xi \right) $ in $E,$ for all $\xi
\in R^{n}$. Hence from $\left( 5.6\right) $ we obtain that the solution of $%
\left( 5.4\right) $ can be represented in the form 
\begin{equation*}
u\left( x\right) =F^{-1}\left( \Psi ^{-1}\left( \xi \right) \hat{f}\left(
\xi \right) \right) .
\end{equation*}%
Moreover, we have 
\begin{equation}
\left\Vert D_{k}^{i}u\right\Vert _{X}=\left\Vert F^{-1}\left( i\xi
_{k}\right) ^{i}\hat{u}\right\Vert _{_{X}}=\left\Vert F^{-1}\xi _{k}^{i}\Psi
^{-1}\left( \xi \right) \hat{f}\right\Vert _{_{X}}  \tag{5.7}
\end{equation}%
and \ 
\begin{equation}
\left\Vert Au\right\Vert _{_{X}}=\left\Vert F^{-1}A\hat{u}\right\Vert
_{_{X}}=\left\Vert F^{-1}A\Psi ^{-1}\left( \xi \right) \hat{f}\right\Vert
_{_{X}}.  \tag{5.8}
\end{equation}

By virtue of $\left( 5.7\right) $ and $\left( 5.8\right) $ for proving $%
\left( 5.5\right) $ it is suffices to show the following estimate%
\begin{equation}
\sum\limits_{k=1}^{n}\sum\limits_{i=0}^{2l_{k}}t_{k}^{\frac{i}{2l_{k}}%
}\left\vert \lambda \right\vert ^{1-\frac{i}{2l_{k}}}\left\Vert \left( i\xi
_{k}\right) ^{i}\hat{u}\right\Vert _{X}+\left\Vert Au\right\Vert _{X}\leq
C\left\Vert f\right\Vert _{X}  \tag{5.9}
\end{equation}%
for all $u\in Y.$ For this aim, it sufficient to show that the operator
functions 
\begin{equation*}
\varphi _{\lambda ,t}\left( \xi \right) =\Psi ^{-1}\left( \xi \right) ,\text{
}\varphi _{ki}\left( \xi \right) =\sum\limits_{i=0}^{2l_{k}}t_{k}^{\frac{i}{%
2l_{k}}}\left\vert \lambda \right\vert ^{1-\frac{i}{2l_{k}}}\xi _{k}^{i}\Psi
^{-1}\left( \xi \right)
\end{equation*}%
are multipliers in $X$ uniformly with respect to $t_{k}$ and $\lambda .$
Firstly, show that $\varphi _{\lambda t}\left( \xi \right) =\Psi ^{-1}\left(
\xi \right) $ \ is a multiplier in $X$ uniformly in $\lambda $ and $t_{k}.$
Indeed, for all $\xi \in R^{n}$ and $\lambda \in S\left( \varphi \right) $
we get 
\begin{equation*}
\left\Vert \varphi _{\lambda t}\left( \xi \right) \right\Vert _{L\left(
E\right) }\leq M\left( 1+\left\vert \lambda +\sum\limits_{k=1}^{n}t_{k}\xi
_{k}^{2l_{k}}\right\vert \right) ^{-1}\leq M_{0}.
\end{equation*}

It is clear that%
\begin{equation*}
\ \frac{\partial }{\partial \xi _{k}}\varphi _{\lambda ,t}\left( \xi \right)
=2l_{k}t_{k}\left[ A+\lambda +\sum\limits_{k=1}^{n}t_{k}\xi _{k}^{2l_{k}}%
\right] ^{-2}\xi _{k}^{2l_{k}-1}\ .
\end{equation*}

Hence, 
\begin{equation}
\left\Vert \xi _{k}\frac{\partial }{\partial \xi _{k}}\varphi _{\lambda
,t}\right\Vert _{L\left( E\right) }\leq 2l_{k}t_{k}\xi
_{k}^{2l_{k}}\left\Vert \left[ A+\lambda +\sum\limits_{k=1}^{n}t_{k}\xi
_{k}^{2l_{k}}\right] ^{-2}\right\Vert \leq  \notag
\end{equation}%
\begin{equation}
2l_{k}t_{k}\xi _{k}^{2l_{k}}\left( 1+\left\vert \lambda
+\sum\limits_{k=1}^{n}t_{k}\xi _{k}^{2l_{k}}\right\vert \right) ^{-2}\leq M.
\tag{5.10}
\end{equation}%
Using the estimate $\left( 5.10\right) $ we show the following uniform
estimate 
\begin{equation}
\left\vert \xi _{1}\right\vert ^{\beta _{1}}...\left\vert \xi
_{n}\right\vert ^{\beta _{n}}\left\Vert D_{\xi }^{\beta }\varphi _{\lambda
,t}\left( \xi \right) \right\Vert _{L\left( E\right) }\leq C  \tag{5.11}
\end{equation}%
for\ $\beta =\beta _{1},...,\beta _{n})\in U_{n}$ \ and $\xi =\left( \xi
_{1},...,\xi _{n}\right) \in V_{n}$. In similar way, we prove that 
\begin{equation}
\left\vert \xi _{1}\right\vert ^{\beta _{1}}...\left\vert \xi
_{n}\right\vert ^{\beta _{n}}\left\Vert D_{\xi }^{\beta }\varphi _{ki}\left(
\xi \right) \right\Vert _{L\left( E\right) }\leq C.  \tag{5.12}
\end{equation}

Since Banach space $E$\ satisfies multiplier condition with respect to $p$
and $\tilde{\gamma},$ in view of estimates $\ \left( 5.11\right) $ and $%
\left( 5.12\right) $ we obtain that the operator-functions $\varphi
_{\lambda ,t},$ $\varphi _{ki}$ are multipliers in\ $X.$ So, we obtain the
estimate $\left( 5.9\right) $ which in turn gives the estimate $\left(
5.5\right) .$ That is we obtain the assertion.

Consider the operator $\tilde{O}_{0}$ in $X$\ generated by the problem $%
\left( 5.4\right) $ that is 
\begin{equation*}
D\left( \tilde{O}_{0}\right) =Y\ \text{and}\ \tilde{O}_{0}u=\sum%
\limits_{k=1}^{n}\left( -1\right) ^{l_{k}}t_{k}D_{k}^{2l_{k}}u+Au.
\end{equation*}

From Theorem 5.1 we obtain

\textbf{Result 5.1. }Assume conditions of Theorem 5.1 are satisfied. Then
the operator $\tilde{O}_{0}$ is positive in $X.$

\textbf{Theorem 5.2.} Suppose the conditions of Theorem 5.1 are satisfied
and $A_{\alpha }\left( x\right) A^{-\left( 1-\left\vert \alpha
:2l\right\vert -\mu \right) }\in L_{\infty }\left( R^{n},L\left( E\right)
\right) $ for $0<\mu <1-\left\vert \alpha \text{:}2l\right\vert $. Then for
all $\ f\in X$ and $\lambda \in S\left( \varphi \right) $\ problem $\left(
5.3\right) $ has a unique solution $u\left( x\right) \in Y$\ and the uniform
coercive estimate holds 
\begin{equation}
\sum\limits_{k=1}^{n}\sum\limits_{i=0}^{2l_{k}}t_{k}^{\frac{i}{2l_{k}}%
}\left\vert \lambda \right\vert ^{1-\frac{i}{2l_{k}}}\left\Vert
D_{k}^{i}u\right\Vert _{X}+\left\Vert Au\right\Vert _{X}\leq C\left\Vert
f\right\Vert _{X}.  \tag{5.13}
\end{equation}%
\ \textbf{Proof. }Consider the problem $\left( 5.3\right) .$ We denote by $%
\tilde{O}_{t}$ the operator in\ $L_{\mathbf{p,\tilde{\gamma}}}\left(
R^{n};E\right) $ generated by problem $\left( 5.3\right) .$ Namely 
\begin{equation}
D\left( \tilde{O}_{t}\right) =Y,\text{ \ }\tilde{O}_{t}u=\tilde{O}_{0}u+%
\tilde{O}_{1}u,  \tag{5.15}
\end{equation}%
where%
\begin{equation*}
\ \tilde{O}_{1}u=\sum\limits_{\left\vert \alpha :2l\right\vert
<1}\prod\limits_{k=1}^{n}t_{k}^{\frac{\alpha _{k}}{2l_{k}}}A_{\alpha
}D^{\alpha }u.
\end{equation*}%
The estimate $\left( 5.5\right) $\ implies that the operator $\tilde{O}%
_{0}+\lambda $ has a bounded inverse from $X$ into $Y$. By Theorem 3.1 for
all $u\in Y$ we get\ 

\begin{equation}
\left\Vert \tilde{O}_{1}u\right\Vert _{X}\leq \sum\limits_{\left\vert \alpha
:2l\right\vert <1}\prod\limits_{k=1}^{n}t_{k}^{\frac{\alpha _{k}}{2l_{k}}%
}\left\Vert A_{\alpha }\left( x\right) D^{\alpha }u\right\Vert _{X}\leq
\sum\limits_{\left\vert \alpha :2l\right\vert
<1}\prod\limits_{k=1}^{n}t_{k}^{\frac{\alpha _{k}}{2l_{k}}}\left\Vert
A^{1-\left\vert \alpha :2l\right\vert -\mu }D^{\alpha }u\right\Vert _{X}\leq
\notag
\end{equation}%
\begin{equation}
C\left[ h^{\mu }\ \left( \sum\limits_{k=1}^{n}t_{k}\left\Vert
D_{k}^{2l_{k}}u\right\Vert _{X}+\left\Vert Au\right\Vert _{X}\right)
+h^{-\left( 1-\mu \right) }\left\Vert u\right\Vert _{X}\right] .  \tag{5.16}
\end{equation}%
Then from $\left( 5.16\right) $ for $u\in Y$ we\ obtain

\begin{equation}
\left\Vert \tilde{O}_{1}u\right\Vert _{X}\leq C\left[ h^{\mu }\left\Vert (%
\tilde{O}_{0}+\lambda )u\right\Vert _{X}+h^{-\left( 1-\mu \right)
}\left\Vert u\right\Vert _{X}.\right]  \tag{5.17}
\end{equation}%
It is clear that 
\begin{equation*}
\ \left\Vert u\right\Vert _{X}=\frac{1}{\lambda }\left\Vert \left( \tilde{O}%
_{0}+\lambda \right) u-\tilde{O}_{0}u\right\Vert _{X},\text{ }u\in Y.
\end{equation*}

\ \ \ By Definition1 we get 
\begin{equation}
\left\Vert u\right\Vert _{X}\leq \frac{1}{\left\vert \lambda \right\vert }%
\left\Vert \left( \tilde{O}_{0}+\lambda \right) u\right\Vert _{X}+\left\Vert 
\tilde{O}_{0}u\right\Vert _{X}\leq \frac{1}{\left\vert \lambda \right\vert }%
\left\Vert \left( \tilde{O}_{0}+\lambda \right) u\right\Vert _{X}+  \notag
\end{equation}%
\begin{equation}
+\frac{1}{\left\vert \lambda \right\vert }\left[ \sum\limits_{k=1}^{n}t_{k}%
\left\Vert D_{k}^{2l_{k}}u\right\Vert _{X}+\left\Vert Au\right\Vert _{X}%
\right] .  \tag{5.18}
\end{equation}%
From Theoem 5.1 and estimates $\left( 5.16\right) -\left( 5.18\right) $ for $%
u\in Y$ we\ obtain 
\begin{equation}
\left\Vert \tilde{O}_{1}u\right\Vert _{X}\leq Ch^{\mu }\left\Vert \left( 
\tilde{O}_{0}+\lambda \right) u\right\Vert _{X}+C_{1}\left\vert \lambda
\right\vert ^{-1}h^{-\left( 1-\mu \right) }\left\Vert \left( \tilde{O}%
_{0}+\lambda \right) u\right\Vert _{X}.  \tag{5.19}
\end{equation}%
Then choosing\ $h$ and $\lambda $ such that \ $Ch^{\mu }<1,$ \ $%
C_{1}\left\vert \lambda \right\vert ^{-1}h^{-\left( 1-\mu \right) }<1$ \
from $\left( 5.19\right) $ we obtain that 
\begin{equation}
\ \ \left\Vert \tilde{O}_{1}\left( \tilde{O}_{0}+\lambda \right)
^{-1}\right\Vert _{L\left( X\right) }<1.  \tag{5.20}
\end{equation}%
Using the relation $\left( 5.15\right) $, Theorem 5.1, $\left( 5.20\right) $
and perturbation theory of linear operators (see for instance $\left[ 8%
\right] $ ) we obtain that the operator $\tilde{O}$\ $+\lambda $ is
invertiable from $X$\ into $Y$. It is implies that for all $f\in X$ problem
\ $\left( 5.3\right) $ have a unique solution $u\in Y$ \ and the estimate $%
\left( 5.13\right) $ holds.

Let $O_{t}$ denotes the operator in $L_{\mathbf{p}}\left( R^{n},E\right) $
generated by problem $\left( 5.1\right) $, i.e. 
\begin{equation*}
D\left( O_{t}\right) =W_{\mathbf{p,\gamma }}^{\left[ 2l\right] }\left(
R^{n},E\left( A\right) ,E\right) ,
\end{equation*}

\begin{equation*}
\ O_{t}u=\sum\limits_{k=1}^{n}(-1)^{l_{k}}t_{k}D_{k}^{\left[ 2l_{k}\right]
}u+Au+\sum\limits_{\left\vert \alpha :2l\right\vert
<1}\prod\limits_{k=1}^{n}t_{k}^{\frac{\alpha _{k}}{2l_{k}}}A_{\alpha }D^{%
\left[ \alpha \right] }u.
\end{equation*}

From Theorem 5.1 and Remark 5.1 we obtain the following

\textbf{Result 5.2. }Assume conditions of Theorem 5.2 are satisfied. Then
for all $f\in L_{\mathbf{p}}\left( R^{n};E\right) $ and $\lambda \in S\left(
\varphi \right) $ the\ equation $\left( 5.1\right) $ has a unique solution $%
u\left( x\right) $ that belongs to $W_{\mathbf{p},\gamma }^{\left[ 2l\right]
}\left( R^{n},E\left( A\right) ,E\right) .$\ Moreover, the uniform coercive
estimate holds 
\begin{equation}
\sum\limits_{k=1}^{n}\sum\limits_{i=0}^{2l_{k}}t_{k}^{\frac{i}{2l_{k}}%
}\left\vert \lambda \right\vert ^{1-\frac{i}{2l_{k}}}\left\Vert D_{k}^{\left[
i\right] }u\right\Vert _{L_{\mathbf{p}}\left( R^{n};E\right) }+\left\Vert
Au\right\Vert _{L_{\mathbf{p}}\left( R^{n};E\right) }\leq C\left\Vert
f\right\Vert _{L_{\mathbf{p}}\left( R^{n};E\right) }.  \tag{5.21}
\end{equation}

\textbf{Result 5.3. }Assume the conditions of Theorem 5.2 are satisfied.
Then the resolvent of operator $O_{t}$ satisfies the following sharp uniform
coercive estimate 
\begin{equation*}
\sum\limits_{k=1}^{n}\sum\limits_{i=0}^{2l_{k}}t_{k}^{\frac{i}{2l_{k}}%
}\left\vert \lambda \right\vert ^{1-\frac{i}{2l_{k}}}\left\Vert D_{k}^{\left[
i\right] }\left( O_{t}+\lambda \right) ^{-1}\right\Vert _{B_{p}}+\left\Vert
A\left( O_{t}+\lambda \right) ^{-1}\right\Vert _{B_{p}}\leq C,
\end{equation*}%
where 
\begin{equation*}
B_{p}=L\left( L_{\mathbf{p}}\left( R^{n};E\right) \right) .
\end{equation*}

From the Result 5.3 and theory of semiroup (see e.g. $\left[ \text{22, \S %
1.14.5}\right] $) we obtain

\textbf{Result 5.4. }Assume conditions of Theorem 5.2 are satisfied for $%
\varphi \in \left( \frac{\pi }{2},\pi \right) .$ Then the operator $O_{t}$
is a generator of analytic semigroup in $X$.

\ \ \ \ \ \textbf{Remark 5.2}. There are a lot of positive operators in the
different concrete Banach spaces. Therefore, putting the concrete Banach
spaces instead of\ $E$ and concrete positive differential, psedodifferential
operators, or finite, infinite matrices instead of $A$ in $\left( 1.1\right)
,$ by virtue of Theorem 5.2 we obtain the separability properties of
different class of degenerate partial differential equations or system of
equations.

\begin{center}
\textbf{6. Abstract Cauchy problem for anisotropic parabolic equation with
parameters}
\end{center}

Consider now, the Cauchy problem $\left( 1.2\right) .$ In this section, we
obta\i n the existence and uniqueness of the maximal regular solution of
problem $\left( 1.2\right) $ in mixed $L_{\mathbf{p,\tilde{\gamma}}}\left(
R^{n};E\right) $ norms.

Let $O_{\varepsilon }$ denote differential operator generated by problem $%
\left( 1.1\right) $ for $t_{k}=\varepsilon _{k}$, $A_{\alpha }\left(
x\right) =0$ and $\lambda =0,$ where $\varepsilon =\left( \varepsilon
_{1},\varepsilon _{2},...,\varepsilon _{n}\right) .$ Let 
\begin{equation*}
X=L_{\mathbf{p,\tilde{\gamma}}}\left( R^{n};E\right) .
\end{equation*}

\textbf{Theorem 6.1. }Assume $E$ is Banach space satisfying multiplier
condition with respect to $\mathbf{p}=\left( p_{1},p_{2},...,p_{n}\right) $
and weighted function $\tilde{\gamma}$. Suppose $A$ is a $\varphi $-positive
operator in Banach space $E$ for $0\leq \varphi <\pi $. Then the operator $%
O_{\varepsilon }$ is uniformly $R$-positive in $X.$

\textbf{Proof. }The Result 5.2 \ implies that the operator $O_{\varepsilon }$
is uniformly positive in $X$. We have to prove the $R$-boundedness of the
set 
\begin{equation*}
\sigma _{t}\left( \xi ,\lambda \right) =\left\{ \lambda \left(
O_{\varepsilon }+\lambda \right) ^{-1}:\lambda \in S_{\varphi }\right\} .
\end{equation*}%
From the Theorem 5.1 we have 
\begin{equation*}
\lambda \left( O_{\varepsilon }+\lambda \right) ^{-1}f=F^{-1}\Phi
_{\varepsilon }\left( \xi ,\lambda \right) \hat{f}\text{, }
\end{equation*}%
for $f\in X,$ where 
\begin{equation*}
\Phi _{\varepsilon }\left( \xi ,\lambda \right) =\lambda \left(
A+L_{\varepsilon }\left( \xi \right) +\lambda \right) ^{-1},\text{ }%
L_{\varepsilon }\left( \xi \right) =\sum\limits_{k=1}^{n}\varepsilon _{k}\xi
_{k}^{2l_{k}}.
\end{equation*}%
By definition of $R$-boundedness, it is sufficient to show that the operator
function $\lambda \left[ \hat{A}\left( \xi \right) +L_{\varepsilon }\left(
\xi \right) +\lambda \right] ^{-1}$ (depended on variable $\lambda $ and
parameters $\xi ,$ $\varepsilon $ ) is uniformly bounded multiplier in $X.$
In a similar manner as in Theorem 5.1 one can easily show that $\Phi
_{t}\left( \xi ,\lambda \right) $ is uniformly bounded multiplier in $X.$
Then,\ by definition of $R$-boundedness we have 
\begin{equation*}
\int\limits_{0}^{1}\left\Vert \sum\limits_{j=1}^{m}r_{j}\left( y\right)
\lambda _{j}\left( O_{\varepsilon }+\lambda _{j}\right)
^{-1}f_{j}\right\Vert _{X}dy=\int\limits_{0}^{1}\left\Vert
\sum\limits_{j=1}^{m}r_{j}\left( y\right) F^{-1}\Phi _{\varepsilon }\left(
\xi ,\lambda _{j}\right) \hat{f}_{j}\right\Vert _{X}dy=
\end{equation*}

\begin{equation*}
\int\limits_{0}^{1}\left\Vert F^{-1}\sum\limits_{j=1}^{m}r_{j}\left(
y\right) \Phi _{\varepsilon }\left( \xi ,\lambda _{j}\right) \hat{f}%
_{j}\right\Vert _{X}dy\leq C\int\limits_{0}^{1}\left\Vert
\sum\limits_{j=1}^{m}r_{j}\left( y\right) f_{j}\right\Vert _{X}dy
\end{equation*}%
for all $\xi _{1},\xi _{2},...,\xi _{m}\in R^{n}$, $\lambda _{1},\lambda
_{2},...,\lambda _{m}\in S_{\varphi },$ $f_{1,}f_{2},...,f_{m}\in X$, $m\in
N $, where $\left\{ r_{j}\right\} $ is a sequence of independent symmetric $%
\left\{ -1,1\right\} $-valued random variables on $\left[ 0,1\right] $.
Hence, the set $\sigma _{\varepsilon }\left( \xi ,\lambda \right) $ is
uniformly $R$-bounded.

Now we are ready to state the main result of this section. Let $G=\left(
0,T\right) \times R^{n}$ and $\mathbf{\tilde{p}}=\left( \mathbf{p,}%
p_{0}\right) .$ Let 
\begin{equation*}
L_{\mathbf{\tilde{p},\gamma }}\left( G;E\right) =L_{p_{0}}\left(
0,T;X\right) \text{.}
\end{equation*}

\textbf{Theorem 6.2. }Assume the conditions of Theorem 6.1 are satisfied for 
$\varphi \in \left( \frac{\pi }{2},\pi \right) $. Then for $f\in L_{\mathbf{%
\tilde{p},\gamma }}\left( G;E\right) $\ problem $\left( 1.2\right) $ has a
unique solution $u\in $ $W_{\mathbf{\tilde{p},\gamma }}^{1,2l}\left(
G;E\left( A\right) ,E\right) $ and the following uniform coercive estimate
holds 
\begin{equation*}
\left\Vert \frac{\partial u}{\partial t}\right\Vert _{L_{\mathbf{\tilde{p}%
,\gamma }}\left( G;E\right) }+\sum\limits_{k=1}^{n}\varepsilon
_{k}\left\Vert D_{k}^{2l_{k}}u\right\Vert _{L_{\mathbf{\tilde{p},\gamma }%
}\left( G;E\right) }+\left\Vert Au\right\Vert _{L_{\mathbf{\tilde{p},\gamma }%
}\left( G;E\right) }\leq C\left\Vert f\right\Vert _{L_{\mathbf{\tilde{p}%
,\gamma }}\left( G;E\right) }.
\end{equation*}%
\textbf{Proof. }The problem $\left( 1.2\right) $ can be expressed as the
following Cauchy problem 
\begin{equation}
\frac{du}{dt}+O_{\varepsilon }u\left( t\right) =f\left( t\right) ,\text{ }%
u\left( 0\right) =0.  \tag{6.1}
\end{equation}

Theorem 6.1 implies that the operator $O_{\varepsilon }$ is $R$-positive and
by Result 5.4 it is a generator of an analytic semigroup in $X.$ Then by
virtue of $\left[ \text{24, Theorem 4.2}\right] $ we obtain that for $f\in
L_{p_{0}}\left( 0,T;X\right) $ problem $\left( 6.1\right) $ has a unique
solution $u\in $ $W_{p_{0}}^{1}\left( 0,T;D\left( O_{\varepsilon }\right)
,X\right) $ and the following uniform estimate holds 
\begin{equation}
\left\Vert \frac{du}{dt}\right\Vert _{L_{p_{0}}\left( 0,T;X\right)
}+\left\Vert O_{\varepsilon }u\right\Vert _{L_{p_{0}}\left( 0,T;X\right)
}\leq C\left\Vert f\right\Vert _{L_{p_{0}}\left( 0,T;X\right) }.  \tag{6.2}
\end{equation}

Since $L_{p_{0}}\left( 0,T;X\right) =L_{\mathbf{\tilde{p},\gamma }}\left(
G;E\right) ,$ by Theorem 5.1 we have 
\begin{equation*}
\left\Vert O_{\varepsilon }u\right\Vert _{L_{p_{0}}\left( 0,T;X\right)
}=\left\Vert O_{\varepsilon }u\right\Vert _{L_{\mathbf{\tilde{p}}}\left(
G;E\right) }.
\end{equation*}%
This relation and the estimate $\left( 6.2\right) $ implies the assertion.

Consider now, the Cauchy problem for degenerate parabolic equation 
\begin{equation}
\frac{\partial u}{\partial t}\ +\sum\limits_{k=1}^{n}(-1)^{l_{k}}\varepsilon
_{k}D_{k}^{^{\left[ 2l_{k}\right] }}u+Au=f\left( t,x\right) ,\text{ } 
\tag{6.3}
\end{equation}

\begin{equation*}
u\left( 0,x\right) =0\text{, }x\in R^{n},\text{ }t\in \left( 0,T\right) ,
\end{equation*}%
where $A$ is a linear operator in a Banach space $E$ and $\varepsilon _{k}$
are small positive parameters.

From Theorem 6.2, and Remark 5.1 we obtain

\textbf{Result 6.1. }Assume conditions of Theorem 6.1 are satisfied for $%
\varphi \in \left( \frac{\pi }{2},\pi \right) $. Then for $f\in L_{\mathbf{%
\tilde{p}}}\left( G;E\right) $\ problem $\left( 6.3\right) $ has a unique
solution $u\in $ $W_{\mathbf{p,\gamma }}^{1,\left[ 2l\right] }\left(
G;E\left( A\right) ,E\right) $ and the following uniform coercive estimate
holds 
\begin{equation*}
\left\Vert \frac{\partial u}{\partial t}\right\Vert _{L_{\mathbf{\tilde{p}}%
}\left( G;E\right) }+\sum\limits_{k=1}^{n}\varepsilon _{k}\left\Vert D_{k}^{%
\left[ 2l_{k}\right] }u\right\Vert _{L_{\mathbf{\tilde{p}}}\left( G;E\right)
}+\left\Vert Au\right\Vert _{L_{\mathbf{\tilde{p}}}\left( G;E\right) }\leq
C\left\Vert f\right\Vert _{L_{\mathbf{\tilde{p}}}\left( G;E\right) }.
\end{equation*}

\begin{center}
\textbf{7. System of parabolic equation of infinite order with small
parameters}
\end{center}

Consider the infinity systems of Cauchy problem for the degenerate
anisotropic parabolic equation

\begin{equation}
\frac{\partial u_{m}}{\partial t}+\sum\limits_{k=1}^{n}\left( -1\right)
^{l_{k}}\varepsilon _{k}D_{k}^{\left[ 2l_{k}\right] }u_{m}+\sum%
\limits_{j=1}^{\infty }d_{j}\left( x\right) u_{j}=f_{m}\left( t,x\right) 
\text{, }  \tag{7.1}
\end{equation}

\begin{equation}
\text{ }u_{m}\left( 0,x\right) =0,\text{ }x\in R^{n},\text{ }t\in \left(
0,T\right) ,\text{ }m=1,2,...,\infty ,  \tag{7.2}
\end{equation}%
\ \ \ where $d_{j}$ are complex valued functions, $\varepsilon _{k}$ are
small positive parameters

\begin{equation*}
B=\left\{ d_{m}\right\} ,\text{ }d_{m}>0,\text{ }u=\left\{ u_{m}\right\} ,%
\text{ }Bu=\left\{ d_{m}u_{m}\right\} ,\text{ }m=1,2,...\infty ,
\end{equation*}%
\begin{equation*}
\text{ }l_{q}\left( B\right) =\left\{ u\in l_{q},\left\Vert u\right\Vert
_{l_{q}\left( D\right) }=\left\Vert Bu\right\Vert _{l_{q}}=\left(
\sum\limits_{m=1}^{\infty }\left\vert d_{m}u_{m}\right\vert ^{q}\right) ^{%
\frac{1}{q}}<\infty \right\} .
\end{equation*}

In this section we show the following result:

\textbf{Theorem 7.1. }For $f\left( t,x\right) =\left\{ f_{m}\left(
t,x\right) \right\} _{1}^{\infty }\in L_{p}\left( G;l_{q}\right) $ problem $%
\left( 7.1\right) -\left( 7.2\right) $ has a unique solution $u\in $ $W_{%
\mathbf{\tilde{p},}\gamma }^{1,\left[ 2l\right] }\left( G,l_{q}\left(
D\right) ,l_{q}\right) $ and the following coercive uniform estimate holds 
\begin{equation}
\left\Vert \frac{\partial u}{\partial t}\right\Vert _{L_{\mathbf{\tilde{p}}%
}\left( G;l_{q}\right) }+\sum\limits_{k=1}^{n}\varepsilon _{k}\left\Vert
D_{k}^{\left[ 2l_{k}\right] }u\right\Vert _{L_{\mathbf{\tilde{p}}}\left(
G;l_{q}\right) }+\left\Vert Au\right\Vert _{L_{\mathbf{\tilde{p}}}\left(
G;l_{q}\right) }\leq C\left\Vert f\right\Vert _{L_{\mathbf{\tilde{p}}}\left(
G;l_{q}\right) }.  \tag{7.3}
\end{equation}

\ \textbf{Proof.} Assume $E=l_{q}$ and $A$ is such that

\begin{equation*}
A=\left[ d_{m}\delta _{mj}\right] ,\text{ }m\text{, }j=1,2,...\infty .
\end{equation*}%
It is clear that the operator $A$ is $R$-positive in $l_{q}$. Then, from
Result 6. 1 we obtain the assertion. Now, consider the following Cauchy
problem

\bigskip 
\begin{equation}
\frac{\partial u_{m}}{\partial t}+\sum\limits_{k=1}^{n}\left( -1\right)
^{l_{k}}\varepsilon _{k}\frac{\partial ^{2l_{k}}u_{m}}{\partial
x_{k}^{2l_{k}}}+\sum\limits_{j=1}^{\infty }d_{j}\left( x\right)
u_{j}=f_{m}\left( t,x\right) \text{, }  \tag{7.4}
\end{equation}

\begin{equation}
\text{ }u_{m}\left( 0,x\right) =0,\text{ }x\in R^{n},\text{ }%
m=1,2,...,\infty .  \tag{7.5}
\end{equation}

From Theorem 7.1 and Remark 5.1 we obtain

\textbf{Result 7.1. }For $f\left( t,x\right) =\left\{ f_{m}\left( t,x\right)
\right\} _{1}^{\infty }\in L_{p}\left( G;l_{q}\right) $ problem $\left(
7.4\right) -\left( 7.5\right) $ has a unique solution $u\in $ $W_{\mathbf{%
\tilde{p},\gamma }}^{1,2l}\left( G,l_{q}\left( D\right) ,l_{q}\right) $ and
the following coercive uniform estimate holds 
\begin{equation*}
\left\Vert \frac{\partial u}{\partial t}\right\Vert _{L_{\mathbf{\tilde{p}%
,\gamma }}\left( G;l_{q}\right) }+\sum\limits_{k=1}^{n}\varepsilon
_{k}\left\Vert \frac{\partial ^{2l_{k}}u}{\partial x_{k}^{2l_{k}}}%
\right\Vert _{L_{\mathbf{\tilde{p},\gamma }}\left( G;l_{q}\right)
}+\left\Vert Au\right\Vert _{L_{\mathbf{\tilde{p},\gamma }}\left(
G;l_{q}\right) }\leq C\left\Vert f\right\Vert _{L_{\mathbf{\tilde{p},\gamma }%
}\left( G;l_{q}\right) }.
\end{equation*}

\begin{center}
\textbf{References}
\end{center}

\begin{enumerate}
\item H. Amann, Linear and quasi-linear equations 1, Birkhauser (1995).

\item R. S. Adams, Sobolev Spaces, Academic Press, New York (1975).

\item D. L. Burkholder, A geometrical condition that implies the existence
certain singular integral of Banach space-valued functions, Proc. conf.
Harmonic analysis in honor of Antony Zygmuhd, Chicago,1981, Wads Worth,
Belmont, (1983), 270-286.

\item O. V. Besov,V. P. Ilin, S. M. Nikolskii, Integral representations of
functions and embedding theorems, Moscow (1975).

\item D. L. Fernandes, On Fourier multipliers of Banach lattice-valued
functions, Rev.Roum. de Math.Pures et appl. (34) (7) (1989), 635-642.

\item V. S. Gul\i ev, To the theory of multiplicators of Fourier integrals
for Banach-valued functions, Trud. Math. Inst. Steklov, (214),
(1996),157-174.

\item M. Girardi, L. Weis, Operator-valued Fourier multiplier theorems on
Lp(X) and geometry of Banach spaces, Journal of Functional Analysis, 204
(2003), (2), 320-354.

\item T. Kato, Perturbation Theory for Linear Operators, Springer-Verlak,
New York (1966).

\item D. S. Kurtz, R. L. Weeden, Results on weighted norm inequalities for
multipliers, Trans. Amer. Math. Soc., 255,1979, 343-362.

\item P. Kree, Sur les multiplicateurs dans $FL$ aves poids, Annales Ins.
Fourier, Grenoble, 16 (2)(1966), 91-121.

\item V. M. Kokilashvili, S. G. Samko, Operators of harmonic analysis in
weighted spaces with non-standard growth,\ \ Journal of Mathematical
Analysis and Applications, 352 (1))1 (2009), 15--34.

\item J. Lions, J. Peetre, Sur one classe d'espases d'interpolation, IHES
Publ. Math. 19(1964), 5-68.

\item P. I. Lizorkin,V. B. Shakhmurov, Embedding theorems for classes of
vector-valued functions, 1,2, IZV.VUZ. USSR , Math., (1989), 70-78, 47-54, .

\item V. G. Mazya, \ Sobolev Spaces, Springer-Verlag, New York (1985).

\item S. L. Sobolev, Certain applications of functional analysis to
mathematical physics, Novosibirskii (1962).

\item S. L. Sobolev, Embeddig theorems for abstract functions, Dok. Akad.
Nauk USSR, 115(1957), 55-59.

\item H. I. Schmeisser and H.Triebel, Topics in Fourier Analysis and
Function Spaces, Willey, Chichester (1987).

\item V. B. Shakhmurov, Abstract capacity of regions and compact embedding
with applications, Acta Mathematica Scientia, 31(1) (2011), 49-67.

\item V. B. Shakhmurov, Embedding and maximal regular differential operators
in Sobolev-Lions spaces, Acta Mathematica Sinica (English Series), 22(5)
(2006), 1493--1508.

\item V. B. Shakhmurov, Estimates of approximation numbers for embedding
operators and applications, Acta Mathematica Sinica, (English Series), 28(
9)( 2012), 1883-1896.

\item V. B. Shakhmurov, Embedding operators and maximal regular
differential-operator equations in Banach-valued function spaces, Journal of
Inequalities and Applications, 4(2005), 329-345.

\item H. Triebel, Interpolation theory, Function spaces, Differential
operators, North-Holland, Amsterdam (1978).

\item H. Triebel H, Spaces of distributions with weights. Multiplier on $%
L_{p}$ spaces with weights. Math. Nachr\textit{.,} (1977)78, 339-356.

\item L. Weis, Operator-valued Fourier multiplier theorems and maximal $%
L_{p} $ regularity, Math. Ann. 319 (2001), 735-75.

\item S. Y. Yakubov, V. B. Shakhmurov, The embedding theorems in the
anisotropic space of vector-functions, Math. Zamet, 22, (1977), 297-30.

\item S. Yakubov and Ya. Yakubov, Differential-operator Equations. Ordinary
and Partial \ Differential Equations, Chapman and Hall /CRC, Boca Raton,
(2000).
\end{enumerate}

\end{document}